\documentclass[10pt,twoside,a4paper,abstract=on]{scrartcl}

\usepackage{authblk}

\usepackage{geometry}
\usepackage{multicol}

\usepackage[utf8]{inputenc}
\usepackage{hyperref}
\usepackage{caption}

\usepackage{todonotes}

\usepackage[toc,title,page]{appendix}

\usepackage{amsmath,amsfonts,amssymb,amsthm}
\usepackage{bbm}
\usepackage{stackengine} 

\usepackage{enumitem}

\usepackage{tikz}  

\usepackage{algorithm}
\usepackage{algpseudocode}

\newtheorem{Theorem}{Theorem}[section]

\newtheorem{Definition}[Theorem]{Definition}

\newtheorem{Remark}[Theorem]{Remark}

\usepackage[a-1b]{pdfx}

\usepackage{graphicx}
\usepackage{subcaption}

\usepackage[maxbibnames=50]{biblatex}
\definecolor{electricultramarine}{rgb}{0.25, 0.0, 1.0}
\definecolor{ikb}{rgb}{0.0, 0.18, 0.65}
\definecolor{green(colorwheel)(x11green)}{rgb}{0.16, 0.5, 0.0}

\usepackage[normalem]{ulem}
\newcommand*{\fzcst}[1]{\relax\ifmmode\text{\textcolor{green(colorwheel)(x11green)}{\sout{\ensuremath{#1}}}}\else\textcolor{green(colorwheel)(x11green)}{\sout{#1}}\fi}

\usepackage{todonotes,varwidth}

\newcommand*{\ykcst}[1]{\relax\ifmmode\text{\textcolor{ikb}{\sout{\ensuremath{#1}}}}\else\textcolor{ikb}{\sout{#1}}
\fi}

\title{\vspace{-0.2cm}Further Approaches of \\
	Dynamical Low-Rank Approximation for SDEs}
\author[1]{Yoshihito Kazashi}
\author[2]{Fabio Nobile}
\author[2]{Fabio Zoccolan}
\affil[1]{Department of Mathematics,
	The University of Manchester, M13 9PL, UK. email: y.kazashi@manchester.ac.uk}
\affil[2]{Institut de Math\'ematiques, \'Ecole Polytechnique F\'ed\'erale de Lausanne, 1015 Lausanne, Switzerland. email: fabio.nobile@epfl.ch, fabio.zoccolan@epfl.ch}
\date{\today}

\begin{document}
   
   \maketitle
	
	\begin{abstract}
		
		In this article, we propose two other DLRA-type dynamics for stochastic differential equations (SDEs) than the one studied in \cite{kazashi2025dynamical}, derived from a minimization of functionals and (informally) from a Stratonovich formulation, respectively. The former approach resembles the DLRA for SDE system proposed in \cite{cao2018stochastic}. Providing the differentiability of the diffusion, the latter procedure registers an additional term in the drift. Indeed, its derivation exploits the Stratonovich formulation to write stochastic processes on manifold, and, hence, possesses a term that depends on the geometry of the manifold itself. These developments open the debate on which formalism is more suitable and what DLRA for SDEs really is.
			\end{abstract}	
    

	\section*{Introduction}
	\addcontentsline{toc}{section}{Introduction}
	
	The \textit{Dynamical Low-Rank Approximation} (DLRA) is a model order reduction technique that is characterized by a linear combination of bases that span in low-dimensional subspaces, are all time-dependent and do not rely on the true solution of the equation to approximate. These features make DLRA very appealing, as it allows to compute this surrogate completely on-the-fly at a cheap cost.
	
	This formalism was first proposed for matrix ODEs in \cite{koch2007dynamical} and then, thanks to its benefits, it was extensively applied in other contexts with outstanding results, too, for instance random and deterministic PDEs (see e.g.\ \cite{bachmayr2021existence,ceruti2022unconventional,einkemmer2019quasi,einkemmer2021asymptotic,kazashi2021existence,kazashi2021stability,sapsis2009dynamically}). In these contexts, one exploits the differentiability of solutions to derive these equations. A first proposal of DLRA-type relations for stochastic differential equations (SDEs) was proposed in \cite{sapsis2009dynamically}, where dynamically-orthogonal formalism was employed to (informally) derive such relations. More precisely, the solution $X$ is represented as a linear combination of $k$ terms for any random realization $\omega$, namely
	\begin{equation}\label{eq: DO}
		X(t,\omega) = \sum_{i = 1}^{k} U^{i}(t)Y^{i}(t,\omega), \quad t \geq 0,
	\end{equation}
	where $\{U^{i}\}_{i=1,\dots, k}$ denotes the time-dependent deterministic basis, orthonormal with respect to a prescribed inner product, and $\{Y^{i}\}_{i=1,\dots, k}$ designates the associated stochastic basis.
	
	A rigorous derivation and well-posedness analysis of DLRA for SDEs was proposed in \cite{kazashi2025dynamical} under standard conditions of Lipschitzianity and linear-growth bound for the coefficients of the SDEs, and for weaker assumptions in \cite{kazashi2026further}. Moreover, a study of their long-time behavior can be found in \cite{bao2026exponential}. This derivation was based on a Itô's formula argument, as SDEs do not possess time-differentiability. To the best of our knowledge, there exists another formalism, different from the one in \cite{kazashi2025dynamical}, proposed in \cite{cao2018stochastic}, where a measure approach was considered. Furthermore, numerical analysis of time and stochastic discretization of these DLRA-for-SDE settings was presented in \cite{kazashi2026dynamicalpartI,kazashi2026dynamicalpartII}. 
	
	The study and implementation of DLRA for SDEs has been increasingly essential in practical applications. Indeed, the growing complexity of contemporary stochastic models, driven by high-dimensional simulations and real-time data-rich applications has made full-order computations of SDE models increasingly impractical. Time-dependent reduced-order modeling has therefore become essential for enabling efficient simulations while preserving the key stochastic features of the underlying system and adapting instantaneously the dynamics to dramatic changes led by the diffusion noise. In this regard, it is essential to understand if there exist other DLRA formulations and which is the most suitable for practical purposes.
	
	The aim of this paper is to propose additional DLRA-type dynamics derived differently from \cite{kazashi2025dynamical} and \cite{cao2018stochastic}, opening the practical, but also philosophical, debate on which formalism is more appropriate. 
	We pursue the following paths:
	\begin{enumerate}
			\item a new DLRA approach inherent to the minimization of the Gramian of the original SDE is proposed. Similarly, the same idea can be translated into the minimization of the mean and the covariance;
		\item another DLRA approach that exploits the stochastic analysis on manifolds in the Stratonovich formulation for SDEs is presented, in contrast to \cite{kazashi2025dynamical} where the Itô framework is directly employed.
	\end{enumerate}
	
	The results of the first point are summarized in Section \ref{sec: min DLRA}. We start by pointing out that the derivation of DLRA for SDEs was not obtained in the same fashion of the traditional well-posed framework \cite{kazashi2025dynamical}. This 2-component surrogate was built using Itô formula and a consistency argument, assuming that the coefficients of the reduced order model have to match the ones of the general SDE if the latter was low-rank. This treatment does not follow the traditional one of \cite{koch2007dynamical}, as the solution of an SDE is usually not time-differentiable. Indeed, 
	therein the derivation of DLRA concerns finding the surrogate of a chosen rank $k$, which minimizes the residual of a certain differential equation, i.e.\ namely the time derivative of this approximation minus the full force term computed in the surrogate itself. Therefore, to obtain the components that define the DLRA, there is no need to employ the true solution of the approximated process, and, hence, the time-dependent bases of DLRA can be computed on the fly. Moreover, the traditional approach turns out to have an equivalent geometrical formulation: the DLRA is the element of a low-rank manifold such that its time derivative is given by projecting the forcing term in the tangent space of that element. The latter interpretation cannot be derived by the formalism of \cite{kazashi2025dynamical} either. 
	
	Here we obtain DLRA surrogates of the type \eqref{eq: DO} by trying to follow the two interpretations of the traditional approach. The former strategy consists in minimizing the differences between two Gramians: the one obtained by the time-discretized DLRA, considered as an unknown, and the one obtained by an Euler-Maruyama method with the DLRA at the previous point of the time mesh as starting point. Then, the limit for the time-mesh going to zero makes us recover the new DO equations. In contrast to \cite{kazashi2025dynamical}, the equation for $U$ shows an additional term depending on the diffusion. The latter technique shows the same DO equations proposed in \cite{cao2018stochastic}. The same treatment can be pursued by simultaneously minimizing mean and covariance of the two discretizations. We provide these DLRA equations in Sections \ref{sec: continuous time limit}-\ref{sec: mean and cov DLRA}.
	
	The second strategy which is illustrated in Section \ref{sec: strat dlra}, is based on exploiting the Stratonovich calculus \cite{protter2012stochastic}, as time-derivatives of SDEs are admitted in this formulation, to seek the geometrical interpretation of traditional DLRA. New equations are obtained by considering a particle-system matrix SDE whose columns are intended as realizations of a continuous random vector DLRA for SDEs and then by taking the limit for the column dimension going to infinity. We illustrate these DLRA equations in Section \ref{sec: DLRA eq strat}.
	
	Well-posedness arguments and additional extensions for these two new methodologies can still be derived similarly to the results stated in \cite{kazashi2025dynamical,kazashi2026further}.
	Moreover, it turns out that all the proposed approaches satisfy the consistency argument developed in \cite{kazashi2025dynamical}, leaving open the practical question of which methodology is the best and the philosophical dilemma of understanding which approach is correct.
	 
	\section{The DLRA-for-SDE framework of \cite{kazashi2025dynamical}}\label{sec:well-posedness SDE}
	We briefly recall the setting where the well-posedness of the dynamically orthogonal framework for SDE was rigorously analyzed.

	Let us consider a stochastic basis $\left( \Omega, \mathcal{F}, \mathbb{P}, (\mathcal{F}_t)_{t \geq 0} \right)$, where $\Omega$ is the probability domain, $\mathcal{F}$ is a $\sigma$-algebra on $\Omega$, $\mathbb{P}$ is a measure of probability on $\Omega$ and $(\mathcal{F}_t)_{t \geq 0}$ is a standard filtration on the probability space $\left(\Omega, \mathcal{F}, \mathbb{P}\right)$.
	We consider $W$ a real $m$-dimensional $(\mathcal{F}_t)$-Brownian motion, denoted as $W(t) = \left(W_1(t), \ldots , W_m(t) \right)^{\top}.$
	We want to establish suitable Dynamical Low-Rank Approximations (DLRAs) for generic SDEs of the following integral form
	\begin{equation}\label{eq:SDE-int}
		X^{\mathrm{true}}= X^{\mathrm{true}}_0 + \int_{0}^{t}a(s,X^{\mathrm{true}}_s)\mathrm{d}s+\int_{0}^{t}b(s,X^{\mathrm{true}}_s)\mathrm{d}W_s, \quad \forall t \in [0, +\infty),
	\end{equation}
	where $X^{\mathrm{true}}_0 = X^{\mathrm{true}}(0)$, the solution of \eqref{eq:SDE-int} is a vector $X^{\mathrm{true}}(t)=\left(X^{\mathrm{true}}_1(t), \ldots, X^{\mathrm{true}}_d(t)\right)^{\top}$, $d \in \mathbb{N}$, whereas the drift $a\colon[0,\infty)\times\mathbb{R}^{d}\to\mathbb{R}^{d}$ 
	and the diffusion $b\colon[0,\infty)\times\mathbb{R}^{d}\to\mathbb{R}^{d\times m}$
	are measurable between the Borel fields $\left([0,\infty)\times\mathbb{R}^{d};\mathcal{B}([0,\infty)\times\mathbb{R}^{d})\right)$ and $\left(\mathbb{R}^{d};\mathcal{B}(\mathbb{R}^{d})\right)$, $\left(\mathbb{R}^{d\times m};\mathcal{B}(\mathbb{R}^{d \times m})\right)$, where $\mathcal{B}$ denotes the Borel $\sigma$-algebra of a given set.
	
	Our low-rank surrogate approximating \eqref{eq:SDE-int} is defined by the pair $\left(U,Y\right)$, solutions of the so-called DO equations~\cite{kazashi2025dynamical,sapsis2009dynamically}:
	\begin{align}
		{C}_{Y_{t}}\dot{U}_{t} & =\mathbb{E}[Y_{t}a(t,U_{t}^{\top}Y_{t})^{\top}](I_{d\times d}-P_{U_{t}}^{\mathrm{row}}),\label{eq:DLR-eq-U}\\
		\mathrm{d}Y_{t} & =U_{t}a(t,U_{t}^{\top}Y_{t})\,\mathrm{d}t+U_{t}b(t,U_{t}^{\top}Y_{t})\mathrm{d}W_{t},\label{eq:DLR-eq-Y}
	\end{align}
	where $C_{Y_t}:= \mathbb{E}[Y_tY_{t}^{\top}]$ is the Gram matrix (or Gramian) of the stochastic basis $Y_t$, and $P^{\text{row}}_{U_t}$ is the projector matrix onto the vector space $\operatorname{span}\{U^{1}_t, \ldots, U^{k}_t\} \subset \mathbb{R}^{d}$, where $U^{i}_t$ is the $i$-th row of  $U_t$. If $U_t$ has orthonormal rows, then $P^{\text{row}}_{U_t}=U_t^{\top}U_t$ i.e.\ the Gramian linked to $U_t$ has linear independent components. Then, the Dynamical Low-Rank approximation $X$ is defined by the product $X=U^{\top}Y$. Notice that \eqref{eq:DLR-eq-U} and \eqref{eq:DLR-eq-Y} are strongly coupled equations, hence their well-posedness is not trivial. Moreover, the evolution of the deterministic modes in \eqref{eq:DLR-eq-U} depends on the law of the process making \eqref{eq:DLR-eq-U}-\eqref{eq:DLR-eq-Y} of McKean-Vlasov type.
	
	Given a suitable rank-$k$ approximation $X_0 :=U_0^{\top}Y_0$ of the initial condition $X_0^{\text{true}}$ of \eqref{eq:SDE-int}, whose rank is assumed to be at least $k$, in \cite{kazashi2025dynamical} the following definition of solution, namely \textit{strong DO solution} is considered.
	
	\begin{Definition}[strong DO solution of rank $k$]\label{def: D0 sol}
		A function $(U,Y) : [0,T] \to \mathbb{R}^{k\times d} \times L^{2}(\Omega;\mathbb{R}^{k})$ is called a \textit{strong DO solution} of rank $k$ for \eqref{eq:SDE-int} if the following conditions are satisfied:
		\begin{enumerate}
			\item {the initial conditions $(U_0,Y_0 )$ are such that $U_0  \in \mathbb{R}^{k\times d}$ is a matrix with orthonormal rows and $Y_0  \in L^{2}(\Omega;\mathbb{R}^{k})$ has linearly independent components};
			\item the curve $t \to U_t \in \mathbb{R}^{k\times d}$ is absolutely continuous on $[0,T]$ and $U_t\dot{U}^{\top}_t = 0 \in \mathbb{R}^{k \times k}$ for a.e. $t \in [0,T]$;
			\item the curve $t \to Y_t(\omega) \in  \mathbb{R}^{k}$ has almost surely continuous paths on $[0,T]$ and is $\mathcal{F}_{t}$-measurable for all $t \in [0,T]$. Moreover,  for any $t \in [0,T]$ the components $Y^{1}_t,\dots,Y^{k}_t$ are linearly independent in $L^2(\Omega)$;
			\item $U$ satisfies equation \eqref{eq:DLR-eq-U} for a.e.~$t \in [0,T]$ and $Y$ is a strong solution of \eqref{eq:DLR-eq-Y} 
			on $[0,T]$.
		\end{enumerate}
	\end{Definition}
	In \cite{kazashi2025dynamical}, it is shown that a unique strong DO solution $(U,Y)$ exists under Lipschitzianity and linear-growth bound conditions on the drift $a$ and the diffusion $b$, which also guarantee well-posedness of the original problem \eqref{eq:SDE-int} (see e.g.\ \cite[Theorem 3.1]{mao2007stochastic} or \cite[Theorem 2.9]{karatzas2012brownian}).

	\section{New Approaches: a minimization derivation}\label{sec: min DLRA}
	
	In this section, we pursue a strategy to derive possible DLRA-type surrogates for SDEs different to the one exploiting the Itô's formula approach in \cite{kazashi2025dynamical}. More specifically, we are considering a minimization procedure, similar to the one proposed in \cite[Section 1]{koch2007dynamical}, where DLRA equations for matrix ordinary differential equations were derived minimizing the difference between the time derivative and the right-hand side of those equations. 
	
	As pointed out in \cite{kazashi2025dynamical}, Itô SDEs do not own time-differentiability property, hence this strategy is not directly doable. Instead, we discretize our SDE under study first, then, we consider a minimization at the discrete level, and, finally, we take the time-step to zero to recover continuous equations. First, we discretize our SDE via a forward Euler-type method and, to approximate it, we consider our DLRA surrogate as linear combination of generic deterministic and stochastic bases satisfying reasonable properties. Then, to derive equations, we minimize the Euclidean error between the Gramian of the forward Euler SDE and the one of the DLRA. Finally, we will formally take the limit $\Delta t \to 0$ in the found relations to retrieve continuous time-differential equations for $U$ and $Y$. We will divide this section in the aforementioned three steps. These new sought DO relations will diversify from the ones in \eqref{eq:DLR-eq-U}-\eqref{eq:DLR-eq-Y} for the spatial basis, which will evolve due to the diffusion matrix, too. This new derivation links to the DO equations proposed in \cite{cao2018stochastic}.
	
	In addition, the \textit{DLR Projector Splitting for Euler-Maruyama} presented in \cite{kazashi2026dynamicalpartI} will result to be an admissible discretization of these new equations. As its numerically performance seems to be superior to other DLR algorithms, it goes without saying that the study of this new approach is not only interesting from the theoretical analysis point of view, but also for applications.
	
	From now on, we suppose $X_n^{\mathrm{DLRA}}$ is the DLRA discrete solution computed at a certain point of the time-mesh $t_n$ with step size $\Delta t_n$. Furthermore, we assume that the DLRA $X_n^{\mathrm{DLRA}}$ can be obtained by the product of discretization of deterministic and stochastic basis $U_n$ and $Y_n$, respectively, with the usual properties, namely $X_n^{\mathrm{DLRA}} = U_n^{\top}Y_n$. 
	
	\subsection{Discretization}
	
	Here, a discretization perspective is considered, which we briefly summarize right away before proceeding in a more detailed treatment. We determine the equations for the bases $\{U_{n}\}_n$ and $\{Y_{n}\}_n$ in the following fashion. First we define $X_{n+1}^{\mathrm{DLRA}} = U_{n+1}^{\top}Y_{n+1}$  at time $t_{n+1}$. Then, we compute the element $X_{n+1}$, \emph{the full-order approximation}, obtained by a first order approximation in time of the true solution at time $t_{n+1}$ computed from the point $X_n=X_n^{\mathrm{DLRA}}$. In detail, $X_{n+1}$ will be defined as the Euler-Maruyama update starting from the point $X_n=X_n^{\mathrm{DLRA}}$.
	Finally, after having preselected a precise norm or distance, we minimize the difference between the Gramians of $X_{n+1}^{\mathrm{DLRA}}$ and $X_{n+1}$ with respect to the possible increments of the deterministic and stochastic basis, $\Delta U_n$ and $\Delta Y_n$, respectively, between $t_n$ and $t_{n+1}$. This procedure will make us recover equations for $U_{n+1}$ and $Y_{n+1}$. In order to obtain differential equations that a candidate continuous DO solution has to satisfy, we will get the limit for the time-step $\Delta t_n$ going to zero of these discrete equations  for $U_n$ and $Y_n$.
	
	In detail, let us consider a partition $\Delta := \left\{t_n \ : \ 0= t_0 < t_1 < \ldots < t_{N-1} < t_{N} = T\right\}$ of $[0,T]$,
	we seek suitable approximations $(U_n)_n$, $(Y_n)_n$ of $(U(t_n))_n$ and $(Y(t_n))_n$, respectively, and define
	an approximate DLRA solution $X_n^{\mathrm{DLRA}} = U_n^{\top}Y_n \approx X(t_n)$. 
	
	\subsubsection{Full-order approximation}\label{sec: full-order}
	The approximation of the true solution $X_{n+1}$ at time $t_{n+1}$ is obtained by a step of the namely Euler-Maruyama method starting with the DLRA solution $X_n^{\mathrm{DLRA}}$ at the time $t_{n}$. This strategy translates into:
	\begin{equation}\label{eq: full-order FE}
		\begin{aligned}
		X_{n+1} &:= X_n^{\mathrm{DLRA}} + a(t_n, X_n^{\mathrm{DLRA}})\Delta t_n + b(t_n, X_n^{\mathrm{DLRA}})\Delta W_n\\
		&= X_n^{\mathrm{DLRA}} + a_n \Delta t_n + b_n \Delta W_n,
		\end{aligned}
	\end{equation}
	where $\Delta t_n = t_{n+1}-t_n$ and $\Delta W_n: = W(t_{n+1})-W(t_{n}) \sim \mathcal{N}(0,\Delta t_n I_{m \times m})$ is the Brownian increment between mesh times $t_n$ and $t_{n+1}$. Notice that by construction $\Delta W_n$ is independent of $X_n^{\mathrm{DLRA}}$. For the sake of notation, hereafter we write $a_n :=a(t_n, X_n^{\mathrm{DLRA}})$ and $b_n := b(t_n, X_n^{\mathrm{DLRA}})$, and, hence, $X_{n+1}  = X_n^{\mathrm{DLRA}} + a_n \Delta t_n + b_n \Delta W_n$.
	
	The Gramian $C_{n+1}$ associated to $X_{n+1}$ defined in \eqref{eq: full-order FE} reads as
	\begin{equation}\label{eq: Gramian X}
		\begin{aligned}
			C_{n+1} =& \mathbb{E}[X_{n+1} X_{n+1}^\top] \\
			= & \mathbb{E}[X_n^{\mathrm{DLRA}} (X_n^{\mathrm{DLRA}})^{\top}] + \mathbb{E}[a_n (X_n^{\mathrm{DLRA}})^{\top}] \Delta t_n \\
			&+ \mathbb{E}[X_n^{\mathrm{DLRA}} a_n^\top] \Delta t_n + \mathbb{E}[a_n a_n^\top] \Delta t_n^2 + \mathbb{E}[b_n b_n^\top] \Delta t_n \\
			=&  C_n^{\mathrm{DLRA}} + \mathbb{E}[a_n (X_n^{\mathrm{DLRA}})^{\top}] \Delta t_n
			+ \mathbb{E}[X_n^{\mathrm{DLRA}} a_n^\top] \Delta t_n + \mathbb{E}[a_n a_n^\top] \Delta t_n^2 + \mathbb{E}[b_n b_n^\top] \Delta t_n,
		\end{aligned}
	\end{equation}
	where in the first line we use the independence of the increment $\Delta W_n$, and $C_n^{\mathrm{DLRA}}$ denotes the Gramian of the DLRA solution $X_n^{\mathrm{DLRA}}$. Notice that $C_{n+1}$ completely characterizes the two moments of $X_{n+1}$.
	
	\subsubsection{Definition of DLRA}
	In this paragraph, we provide a discussion on how we want to define our DLRA surrogate.
	Starting from the DLRA $X_n^{\mathrm{DLRA}}$, with $X_n^{\mathrm{DLRA}} = U_n^{\top}Y_n$, we aim to find equations for $U_{n+1}$ and $Y_{n+1}$ to obtain the DLRA $X_{n+1}^{\mathrm{DLRA}}$, with $X_{n+1}^{\mathrm{DLRA}}= U_{n+1}^{\top}Y_{n+1}$ at the time $t_{n+1}$. We ask for the following discrete updates
    \begin{equation}\label{eq: increments}
    	U_{n+1}:= U_n + \Delta U_n , \quad Y_{n+1} = Y_n + \Delta Y_n,
    \end{equation}
    where $U_n \in \mathbb{R}^{k \times d}$ a deterministic matrix with $U_n U_n^\top = I_{k \times k}$, i.e.\ $U_n$ has orthonormal rows, and $U_n \Delta U_n^\top = 0$, whereas $Y_n \in L^2(\Omega,\mathbb{R}^k)$ and $\Delta Y_n$ is adapted to the given filtration. Notice that for $\Delta t_n \to 0$ we recover the properties $2.,3.$ of Definition \ref{def: D0 sol}. The hypothesis $U_n \Delta U_n^\top = 0$ is assumed in order to recover the continuous gauge condition when $\Delta t_n \to 0$.
	To have a direct connection with the continuous DLRA, we will call the property $U_n \Delta U_n^\top = 0$ as the \emph{discrete gauge condition}.
	
    The Gramian $C_{n+1}^{\mathrm{DLRA}}$ associated to the DLRA solution $X_{n+1}^{\mathrm{DLRA}}$ at time $t_{n+1}$ reads as 
	\begin{equation*}
		C_{n+1}^{\mathrm{DLRA}} = U_{n+1}^\top C_{Y_{n+1}} U_{n+1},
	\end{equation*}
	where $C_{Y_n}= \mathbb{E}[Y_nY_n^{\top}] \in \mathbb{R}^{k \times k}$ is the Gramian of the stochastic basis $Y_n$.
	The Gramian $C_{Y_{n+1}} = \mathbb{E}[Y_{n+1}Y_{n+1}^{\top}]$ of $Y_{n+1}$ can be rewritten as the Gramian at the previous step plus a matrix increment, i.e.\ $C_{Y_n} + \Delta C_{Y_n}$, where $\Delta C_{Y_n} \in \mathbb{R}^{k \times k}$. For the ease of computation, we will derive equations for $\Delta Y_n$ from the one of $\Delta C_{Y_n}$ in a later moment. Therefore, one has
	\begin{equation}\label{eq: C DLRA}
		\begin{aligned}
		C_{n+1}^{\mathrm{DLRA}} &\approx (U_n + \Delta U_n)^\top (C_{Y_n} + \Delta C_{Y_n}) (U_n + \Delta U_n) \\
		& = U_n^\top C_{Y_n} U_n + (\Delta U_n)^\top C_{Y_n} U_n + U_n^\top C_{Y_n} \Delta U_n + U_n^\top \Delta C_{Y_n} U_n,
				\end{aligned}
	\end{equation}
	where here and from now on we consider higher order terms involved in the expression of $C_{n+1}^{\mathrm{DLRA}}$ as negligible with respect to $\Delta t_n$.
	To obtain equations for $U_{n+1}$ and $Y_{n+1}$, we have to determine equations for $\Delta U_n$ and $\Delta Y_n$ (and, hence, $\Delta C_{Y_n}$). In order to find these relations, we minimize the Frobenius distance between the Gramian $C_{n+1}$ defined in \eqref{eq: Gramian X} and $C_{n+1}^{\mathrm{DLRA}}$, similarly to the treatment of \cite[Section 1]{koch2007dynamical}, as detailed in the following subsection.
	\subsection{Minimization of the Frobenius distance between the two Gramians} 
	The minimization of the Frobenius distance between $C_{n+1}$ and $C_{n+1}^{\mathrm{DLRA}}$ with respect to $\Delta U_n$ and $\Delta Y_n$ is described as follows. We seek $(\Delta U_n,\Delta Y_n)$ satisfying
	\begin{equation}\label{eq: min Delta U and Y}
		\min_{\Delta U_n , \Delta Y_n} \|C_{n+1}^{\mathrm{DLRA}}-C_{n+1}\|_{\mathrm{F}}=!.
    \end{equation}	
    Notice that \eqref{eq: min Delta U and Y} is equivalent to
    \begin{equation}\label{eq: min Delta U and Y tr}
    	\min_{\Delta U_n , \Delta Y_n}  \|C_{n+1}^{\mathrm{DLRA}}-C_{n+1}\|_{\mathrm{F}}^2 = \min_{\Delta U_n , \Delta Y_n}  \mathrm{Tr}\left( \left(C_{n+1}^{\mathrm{DLRA}}-C_{n+1} \right) \left(C_{n+1}^{\mathrm{DLRA}}-C_{n+1}\right)^{\top} \right)=!,
    \end{equation}	
    hence, actually solving a first optimality problem, the minimum is reached via setting the gradient of $\mathrm{Tr}\left( \left(C_{n+1}^{\mathrm{DLRA}}-C_{n+1} \right) \left(C_{n+1}^{\mathrm{DLRA}}-C_{n+1}\right)^{\top} \right)$ with respect to $\Delta U_n$ and $\Delta Y_n$ equal to 0 in a variational formulation.
    
	To compute the derivatives that will appear in the variational formulation, the following matrix identities will be useful in future computations \cite{petersen2008matrix}: for all matrices $A,B$ with suitable dimensions
    \begin{equation}\label{eq: matrix der}
    	\begin{aligned}
    	\frac{\partial}{\partial X} \operatorname{Tr}(XA) = A^{\top},\quad
    	\frac{\partial}{\partial X} \operatorname{Tr}(AXB) &= A^{\top} B^{\top}, \quad
    	\frac{\partial}{\partial X} \operatorname{Tr}(AX^{\top}B) = BA  \\
    	\frac{\partial}{\partial X} \operatorname{Tr}(X^{\top}A) &= A, \quad 
    	\frac{\partial}{\partial X} \operatorname{Tr}(AX^{\top}) = A.
    	\end{aligned}
    \end{equation}
    From \eqref{eq: C DLRA}, one has that 
	\begin{equation}\label{eq: Gramian DLRA}
			\begin{aligned}
		C_{n+1}^{\mathrm{DLRA}} (C_{n+1}^{\mathrm{DLRA}})^\top \approx & U_n^\top C_{Y_n} C_{Y_n} U_n^\top + U_n^\top C_{Y_n} C_{Y_n} \Delta U_n + U_n^\top C_{Y_n} \Delta C_{Y_n} U_n \\
			&+ (\Delta U_n)^\top C_{Y_n} C_{Y_n} U_n + \Delta U_n^\top C_{Y_n} C_{Y_n} \Delta U_n + \Delta U_n^\top C_{Y_n} \Delta C_{Y_n} U_n\\
			& + U_n^\top C_{Y_n} \Delta U_n \Delta U_n^\top C_{Y_n} U_n + U_n^\top \Delta C_{Y_n}  C_{Y_n} U_n  + U_n^\top \Delta C_{Y_n} C_{Y_n} \Delta U_n \\
			& + U_n^\top \Delta C_{Y_n} \Delta C_{Y_n} U_n, \\
		\end{aligned}
	\end{equation}
	where we exploited the orthogonality of the rows of $U_n$ and the discrete gauge condition.
	The derivative of \eqref{eq: Gramian DLRA} with respect to the $\Delta U_n$ increment reads as
	\begin{equation}
		\begin{aligned}
		\frac{\partial \operatorname{Tr}\left( C_{n+1}^{\mathrm{DLRA}} (C_{n+1}^{\mathrm{DLRA}})^\top \right) }{\partial \Delta U_n}
		& =  2C_{Y_n} C_{Y_n} U_n + 4C_{Y_n} C_{Y_n} \Delta U_n + 2C_{Y_n} \Delta C_{Y_n} U_n .
		\end{aligned}
	\end{equation}
	
	Furthermore, using the discrete gauge condition and the orthogonality of $U_n$, the derivative of \eqref{eq: Gramian DLRA} with respect to the $\Delta C_{Y_n}$ increment is
	\begin{equation}
		\begin{aligned}
		\frac{\partial \operatorname{Tr} \left(  C_{n+1}^{\mathrm{DLRA}}(C_{n+1}^{\mathrm{DLRA}})^\top \right)}{\partial \Delta C_{Y_n}} 
			&= 2 C_{Y_n} + 2 \Delta C_{Y_n}.
		\end{aligned}
	\end{equation}
	
	Moreover,
		\begin{equation}
		\begin{aligned}
			\frac{\partial \operatorname{Tr}\left( C_{n+1} (C_{n+1}^{\mathrm{DLRA}})^\top \right) }{\partial \Delta U_n}
			& = \frac{\partial \operatorname{Tr}\left( C_{n+1} \left(U_n^\top C_{Y_n} U_n + U_n^{\top} C_{Y_n} \Delta U_n + (\Delta U_n)^\top C_{Y_n} U_n + U_n^\top \Delta C_{Y_n} U_n \right)^{\top} \right)}{\partial \Delta U_n}\\
			&= 2 C_{Y_n} U_n C_{n+1}
		\end{aligned}
	\end{equation}
	and
		\begin{equation}
		\begin{aligned}
			\frac{\partial \operatorname{Tr}\left( C_{n+1} (C_{n+1}^{\mathrm{DLRA}})^\top \right) }{\partial \Delta C_{Y_n}}
			&= U_n C_{n+1} U_n^{\top}.
		\end{aligned}
	\end{equation}
	
	Therefore, first we have
	\begin{equation*}
		\begin{aligned}
			\frac{\partial \operatorname{Tr}\left( \left(C_{n+1} - C_{n+1}^{\mathrm{DLRA}} \right)\left(C_{n+1} - C_{n+1}^{\mathrm{DLRA}} \right)^{\top} \right)}{\partial \Delta C_{Y_n}} &= \frac{\partial \operatorname{Tr} \left( -2 C_{n+1} (C_{n+1}^{\mathrm{DLRA}})^\top + C_{n+1}^{\mathrm{DLRA}} (C_{n+1}^{\mathrm{DLRA}})^\top \right) }{\partial \Delta C_{Y_n}} \\
			&= - 2 U_n C_{n+1} U_n^{\top} + 2 C_{Y_n} + 2 \Delta C_{Y_n}  \\
		\end{aligned}
	\end{equation*}
	Assuming that $C_{n+1}$ is of rank at least $k$ and $C_{Y_n}$, $C_{n+1}$ are of full rank, then $\Delta C_{Y_n}$ is obtaining projecting the full order Gramian $C_{n+1}$ and subtracting the Gramian of the stochastic basis of the DLRA at time $t_n$:
	\begin{equation}\label{eq: Delta C_y}
				\Delta C_{Y_n} = U_n C_{n+1} U_n^{\top} - C_{Y_n}.
	\end{equation}
	
	\subsubsection{The equation for $\Delta U_n$}
	To derive the deterministic increment $\Delta U_n$, similarly we get
	\begin{equation*}
	\begin{aligned}
		\frac{\partial \operatorname{Tr}\left( (C_{n+1} - C_{n+1}^{\mathrm{DLRA}})(C_{n+1} - C_{n+1}^{\mathrm{DLRA}})^{\top} \right)}{\partial \Delta U_n}
		&= \frac{\partial \operatorname{Tr}\left( -2 C_{n+1} (C_{n+1}^{\mathrm{DLRA}})^{\top} + C_{n+1}^{\mathrm{DLRA}} (C_{n+1}^{\mathrm{DLRA}})^{\top} \right) }{\partial \Delta U_n} \\
		&= 2 C_{Y_n} C_{Y_n} U_n + 4 C_{Y_n} C_{Y_n} \Delta U_n + 2 C_{Y_n} \Delta C_{Y_n} U_n \\
		&\quad - 4 C_{Y_n} U_n C_{n+1} \\
		&= 2 C_{Y_n} C_{Y_n} U_n + 4 C_{Y_n} C_{Y_n} \Delta U_n + 2 C_{Y_n} U_n C_{n+1} P_{U_n} \\
		&\quad - 2 C_{Y_n} C_{Y_n} U_n - 4 C_{Y_n} U_n C_{n+1},
	\end{aligned}
\end{equation*}
where in the last line we employ relation \eqref{eq: Delta C_y}.

To find $\Delta U_n$, as we are looking for a rectangular matrix, in a variational formulation one seeks
\begin{equation*}
	\langle \frac{\partial \operatorname{Tr}\left( (C_{n+1} - C_{n+1}^{\mathrm{DLRA}})(C_{n+1} - C_{n+1}^{\mathrm{DLRA}})^{\top} \right)}{\partial \Delta U_n}, V \rangle_{\mathrm{F}} =0, \quad \forall V \in \operatorname{Im}({U_n^{\top}}^\perp),
\end{equation*}
where $\mathrm{Im}(A)$ denotes the image of the operator $A$.
 Therefore, assuming $C_{Y_n}$ to be invertible, and $V \in \text{Im}(P_{U_n}^\perp)$, in a variational formulation one gets that 
	\begin{equation*}
		\begin{aligned}
			0=&\langle 2 C_{Y_n} U_n+ 4 C_{Y_n} \Delta U_n + 2 U_n C_{n+1} P_{U_n} - 4 U_n C_{n+1}, V \rangle_{\mathrm{F}} = 0 \\
			=&\langle 4 C_{Y_n} \Delta U_n - 4 U_n C_{n+1}, V \rangle_{\mathrm{F}} = 0,\quad \forall V \in \mathbb{R}^{d \times k}, \ V \in \operatorname{Im}({U_n^{\top}}^\perp); \\
			\end{aligned}
\end{equation*}
or equivalently 
\begin{equation*}
 4 C_{Y_n} \Delta U_n P_{U_n}^\perp - 4 U_n C_{n+1} P_{U_n}^\perp = 0,
\end{equation*}
which implies
\begin{equation*}
	 \Delta U_n -  C_{Y_n}^{-1} U_n C_{n+1} P_{U_n}^\perp = 0,
\end{equation*}
and, hence, 
	\begin{equation*}
	\begin{aligned}
			 \Delta U_n P_{U_n}^\perp = \Delta U_n = C_{Y_n}^{-1} U_n C_{n+1} P_{U_n}^\perp, 
		\end{aligned}
	\end{equation*}
where we use the fact that $P_{U_n} V^{\top} = 0$ as  $V \in \operatorname{Im}({U_n^{\top}}^\perp)$.
	Considering the relation \eqref{eq: Gramian X}, then the update $\Delta U_n$ is determined by
	\begin{equation*}
		\begin{aligned}
               \Delta U_n &= C_{Y_n}^{-1} U_n C_{n+1} P_{U_n}^\perp\\
                                  &= C_{Y_n}^{-1} U_n \left(C_n^{\mathrm{DLRA}} + \mathbb{E}[a_n (X_n^{\mathrm{DLRA}})^{\top}] \Delta t_n 
                                  + \mathbb{E}[X_n^{\mathrm{DLRA}} a_n^\top] \Delta t_n + \mathbb{E}[a_n a_n^\top] \Delta t_n^2 + \mathbb{E}[b_n b_n^\top] \Delta t_n \right) P_{U_n}^\perp \\
                                  &= C_{Y_n}^{-1} U_n \left(U_n^\top C_{Y_n} U_n+ \mathbb{E}[a_n Y_n^{\top}U_n] \Delta t_n  
                                  + \mathbb{E}[U_n^{\top}Y_n a_n^\top] \Delta t_n + \mathbb{E}[a_n a_n^\top] \Delta t_n^2 + \mathbb{E}[b_n b_n^\top] \Delta t_n \right) P_{U_n}^\perp\\
                                  &= C_{Y_n}^{-1} \mathbb{E}[Y_n a_n^\top] P_{U_n}^\perp \Delta t_n + C_{Y_n}^{-1} U_n \mathbb{E}[a_n a_n^\top] P_{U_n}^\perp \Delta t_n^2 + C_{Y_n}^{-1} U_n \mathbb{E}[b_n b_n^\top] P_{U_n}^\perp \Delta t_n \\
		\end{aligned}
	\end{equation*}
	which, in a first order approximation with respect to $\Delta t_n$, is equivalent to
	\begin{equation}\label{eq: Delta U_n}
		\begin{aligned}
			\Delta U_n & \approx C_{Y_n}^{-1} \mathbb{E}[Y_n a_n^\top] P_{U_n}^\perp \Delta t_n + C_{Y_n}^{-1} U_n \mathbb{E}[b_n b_n^\top] P_{U_n}^\perp \Delta t_n.\\
		\end{aligned}
	\end{equation}

	\subsubsection{The equation for $\Delta Y_n$}
	Now, we want to determine the equation for the increment $\Delta Y_n$ exploiting the relation for $\Delta C_{Y_n}$ obtained in \eqref{eq: Delta C_y}. We have that the Gramian of $Y_{n+1}$ has to satisfy
	\begin{equation*}
		\begin{aligned}
			C_{Y_{n+1}} &= \mathbb{E}[Y_{n+1} Y_{n+1}^{\top}] \\
			&= \mathbb{E}[(Y_n + \Delta Y_n)(Y_n + \Delta Y_n)^{\top}] \\
			&= C_{Y_n} + \mathbb{E}[Y_n \Delta Y_n^{\top}] + \mathbb{E}[\Delta Y_n Y_n^{\top}] + \mathbb{E}[\Delta Y_n \Delta Y_n^{\top}]. \\
		\end{aligned}
	\end{equation*}
	But from \eqref{eq: Delta C_y} one has that
	\begin{equation*}
		 	C_{Y_{n+1}}  = C_{Y_n} + \Delta C_{Y_n} = C_{Y_n} + U_n C_{n+1} U_n^{\top} - C_{Y_n} = U_n C_{n+1} U_n^{\top},
	\end{equation*}
	which using \eqref{eq: Gramian DLRA} and \eqref{eq: increments} implies the following relation
	\begin{equation*}
		\begin{aligned}
	    & \mathbb{E}[Y_n \Delta Y_n^{\top}] + \mathbb{E}[\Delta Y_n Y_n^{\top}] + \mathbb{E}[\Delta Y_n \Delta Y_n^{\top}] \\
			=& U_n C_{n+1} U_n^{\top} - C_{Y_n} \\
			=& C_{Y_n} + U_n \mathbb{E}[a_n Y_n^{\top}] \Delta t_n + \mathbb{E}[Y_n a_n^{\top}] U_n^{\top} \Delta t_n + U_n \mathbb{E}[a_n a_n^{\top}] U_n^{\top} \Delta t_n^2 + U_n \mathbb{E}[b_n b_n^{\top}] U_n^{\top} \Delta t_n - C_{Y_n} \\
			= & U_n \mathbb{E}[a_n Y_n^{\top}] \Delta t_n + \mathbb{E}[Y_n a_n^{\top}] U_n^{\top} \Delta t_n + U_n \mathbb{E}[a_n a_n^{\top}] U_n^{\top} \Delta t_n^2 + U_n \mathbb{E}[b_n b_n^{\top}] U_n^{\top} \Delta t_n.\\
		\end{aligned}
	\end{equation*}
	Moreover, we have that 
		\begin{equation}\label{eq: mean Y_n}
		\begin{aligned}
			\mathbb{E}[Y_{n+1}] = 	\mathbb{E}[Y_{n}]+ \mathbb{E}[\Delta Y_n]  &= (U_n + \Delta U_n) \left(\mathbb{E}[X_n^{\mathrm{DLRA}}] + \mathbb{E}[a_n] \Delta t_n \right) \approx  \mathbb{E}[Y_n] + U_n \mathbb{E}[a_n] \Delta t_n.
		\end{aligned}
	\end{equation}
	in a first order of optimality.
	Matching the conditions \eqref{eq: Delta C_y} and \eqref{eq: mean Y_n}, we find that the equation for $\Delta Y_n$ reads as
	\begin{equation}\label{eq: Delta Y_n}
		\Delta Y_n = U_n a_{n} \Delta t_n + U_n b_n \Delta W_n.
	\end{equation}
	
	\subsection{The DO equation and the relation for $X_n$}
	Therefore, asking for $U_{n+1} = U_n + \Delta U_n$ and $Y_{n+1} = Y_n + \Delta Y_n$ for each $n$, thanks to relations \eqref{eq: Delta U_n} and \eqref{eq: Delta Y_n} we obtained the following discretized DO equations
		\begin{align}
		U_{n+1} &= U_n + C_{Y_n}^{-1} \mathbb{E}[Y_n a_n^\top] P_{U_n}^{\perp} \Delta t_n 
		+ C_{Y_n}^{-1} U_n \mathbb{E}[b_n b_n^\top] P_{U_n}^{\perp} \Delta t_n \label{eq: U_n gramian} \\
		Y_{n+1} &= Y_n + U_n a_n \Delta t_n + U_n b_n \Delta W_n. \label{eq: Y_n gramian}
		\end{align}
	By multiplying $U_{n+1}$ and $Y_{n+1}$ we get
	\begin{equation*}
		\begin{aligned}
			X_{n+1}^{\mathrm{DLRA}} =& X_n^{\mathrm{DLRA}} + P_{U_n} a_n \Delta t_n + P_{U_n} b_n \Delta W_n + P_{U_n}^{\perp} P_{Y_n} a_n \Delta t_n + P_{U_n}^{\perp} \mathbb{E}[a_n Y_n^\top] C_{Y_n}^{-1} U_n a_n \Delta t_n^2 \\
			&+ P_{U_n}^{\perp} \mathbb{E}[a_n Y_n^\top] C_{Y_n}^{-1} U_n b_n \Delta W_n \Delta t_n + P_{U_n}^{\perp} \mathbb{E}[b_n b_n^\top] U_n^\top C_{Y_n}^{-1} Y_n \Delta t_n \\
			&+ P_{U_n}^{\perp} \mathbb{E}[b_n b_n^\top] U_n^\top C_{Y_n}^{-1} U_n a_n \Delta t_n^2 + P_{U_n}^{\perp} \mathbb{E}[b_n b_n^\top] U_n^\top C_{Y_n}^{-1} U_n b_n \Delta W_n \Delta t_n.
		\end{aligned}
	\end{equation*}
	Via discarding all terms of order greater than $O(\Delta t_n)$ (i.e.\ keeping only terms with $\Delta t_n$ and $\Delta W_n$) one obtains
	\begin{equation}\label{eq: new X_{n+1}}
		\begin{aligned}
		X_{n+1}^{\mathrm{DLRA}}	\approx X_n^{\mathrm{DLRA}} &+ P_{U_n} a_n \Delta t_n + P_{U_n} b_n \Delta W_n + P_{U_n}^{\perp} P_{Y_n} a_n \Delta t_n 
			+ P_{U_n}^{\perp} \mathbb{E}[b_n b_n^\top] U_n^\top C_{Y_n}^{-1} Y_n \Delta t_n. 
		\end{aligned}
	\end{equation}
	Notice that \eqref{eq: new X_{n+1}} differs from the approximation obtained by the DLR Projector Splitting for SDE \cite{kazashi2026dynamicalpartI} by the term $P_{U_n}^{\perp} \mathbb{E}[b_n b_n^\top] U_n^\top C_{Y_n}^{-1} Y_n \Delta t_n$. However,  this additional term is still in the tangent space of the manifold of rank $k$ processes at the point $X_n$; indeed it holds that
	\begin{equation*}
		\begin{aligned}
			\left(P_{U_n}^{\perp} P_{Y_n} + P_{U_n}\right)\left[P_{U_n}^{\perp} \mathbb{E}[b_n b_n^\top] U_n^\top C_{Y_n}^{-1} Y_n \Delta t_n\right] 
			& = P_{U_n}^{\perp} P_{Y_n} \left[P_{U_n}^{\perp} \mathbb{E}[b_n b_n^\top] U_n^\top C_{Y_n}^{-1} Y_n \Delta t_n\right] \\
			&= P_{U_n}^{\perp} \mathbb{E}[b_n b_n^\top] U_n^\top C_{Y_n}^{-1} \mathbb{E}[Y_n Y_n^\top] C_{Y_n}^{-1} Y_n \Delta t_n \\
			&= P_{U_n}^{\perp} \mathbb{E}[b_n b_n^\top] U_n^\top C_{Y_n}^{-1} Y_n \Delta t_n.
		\end{aligned}
	\end{equation*}
	 
	 \subsection{Continuous-Time Limit}\label{sec: continuous time limit}
	 Equations \eqref{eq: U_n gramian} and \eqref{eq: Y_n gramian} are the starting point to define new differential equations for DLRA. We suppose that $(U_n)_n$ and $(Y_n)_n$ are sequences approximating a deterministic function $(U_t)_{t\geq 0}$ and a stochastic one $(Y_t)_{t\geq 0}$, respectively, in the time mesh $(t_n)_n$, where $U$,$Y$ satisfy properties $1.,2.,3.$ of Definition \ref{def: D0 sol}. We define the equations that $U$ and $Y$ satisfy as the (informal) limit 
	 for $\Delta t_n \to 0$ of \eqref{eq: U_n gramian} and \eqref{eq: Y_n gramian}. Therefore, for $Y_t$ we obtain that
	 \begin{equation*}
	 	Y_{n+1} = Y_n + U_n \, a_{n} \, \Delta t_n + U_n \, b_n \, \Delta W_n \ \xrightarrow{\Delta t_n \to 0} \ 
	 	\mathrm{d}Y_t = U_t \, a(t,X_t) \, \mathrm{d}t + U_t \, b(t,X_t) \, \mathrm{d}W_t,
	 \end{equation*}
	 whereas for  $ U_t $ one has
	 \begin{equation*}
	 	\begin{aligned}
	 	&U_{n+1} = U_n + C_{Y_n}^{-1} \mathbb{E}[Y_n a_n^{\top}] P_{U_n}^\perp \Delta t_n 
	 	+ C_{Y_n}^{-1} U_n \mathbb{E}[b_n b_n^{\top}] P_{U_n}^\perp \Delta t_n \\
	 	\ \xrightarrow{\Delta t_n \to 0} \ &
	 	\mathrm{d}U_t = C_{Y_t}^{-1} \, \mathbb{E}[Y_t a(t, X_t)^{\top}] \, P_{U_t}^\perp \mathrm{d}t 
	 	+ C_{Y_t}^{-1} U_t \mathbb{E}[b(t, X_t) b(t, X_t)^{\top}] P_{U_t}^\perp \mathrm{d}t.
	 	\end{aligned}
	 \end{equation*}
	 Therefore, the DO system of equations that the pair $(U,Y)$ satisfies is
	 \begin{align}
	 	\mathrm{d}U_t &= C_{Y_t}^{-1} \, \mathbb{E}[Y_t a(t, X_t)^{\top}] \, P_{U_t}^\perp \mathrm{d}t 
	 	+ C_{Y_t}^{-1} U_t \mathbb{E}[b(t, X_t) b(t, X_t)^{\top}] P_{U_t}^\perp \mathrm{d}t \label{eq: U_t gramian}\\
	 	\mathrm{d}Y_t &= U_t \, a(t, X_t)\, \mathrm{d}t + U_t \, b(t, X_t) \, \mathrm{d}W_t \label{eq: Y_t gramian}.
	 \end{align}
	 Notice that \eqref{eq: Y_t gramian} is the same relation of \eqref{eq:DLR-eq-Y}, whereas \eqref{eq: U_t gramian} differs from \eqref{eq:DLR-eq-U} by a term depending on the diffusion component. This implies that the subspace $U_t$ in \eqref{eq: U_t gramian} evolves also giving possible directions from the noise term, unlike \eqref{eq:DLR-eq-U}, where the subspace evolves only considering directions obtained by the drift $a$.
	 
	 Via Itô's formula, we can also derive an equation for the surrogate $X_t$: indeed,
	 \begin{equation}\label{eq: new DLRA min}
	 	\begin{aligned}
	 	\mathrm{d}X_t =& 	\mathrm{d}(U_t^\top Y_t) = (\mathrm{d}U_t^\top) Y_t + U_t^\top \mathrm{d}Y_t\\
	 		= & P_{U_t}^\perp P_{Y_t} [a(t,X_t)]\,\mathrm{d}t + P_{U_t}^\perp \mathbb{E}[(b(t,X_t) b(t,X_t)^\top)] U_t^{\top} C_{Y_t}^{-1} Y_t\,\mathrm{d}t + P_{U_t}a(t,X_t)\,\mathrm{d}t + P_{U_t} b(t,X_t)\,\mathrm{d}W_t\\
	 		= & P_{X_t} [a(t,X_t)]\,\mathrm{d}t + P_{U_t} b(t,X_t)\,\mathrm{d}W_t + P_{U_t} \mathbb{E}[b(t,X_t) b(t,X_t)^\top] U_t^{\top}C_{Y_t}^{-1} Y_t\,dt,
	 	\end{aligned}
	 \end{equation}
	 where $P_{X_t} := P_{U_t}^\perp P_{Y_t} + P_{U_t}$ is the projection onto the tangent of $L^2(\Omega, \mathbb{R}^d)$ processes of rank $k$ in the point $X_t$.
	 
	 One can question the well-posedness of \eqref{eq: U_t gramian}-\eqref{eq: Y_t gramian}. Actually, as \eqref{eq: U_t gramian}-\eqref{eq: Y_t gramian} differs from \eqref{eq:DLR-eq-U}-\eqref{eq:DLR-eq-Y} from a small term, the strategy to prove existence and uniqueness follows closely to discussion developed in \cite{kazashi2025dynamical,kazashi2026further} depending on the given assumptions. For the sake of completeness, we briefly sketch the main points of a possible strategy among the admissible ones:
	 \begin{itemize}
	 	\item finding a time $T>0$ such that the Gram matrices associated to $U$ and $Y$ are always non-degenerate and, hence, equations \eqref{eq: U_t gramian} and \eqref{eq: Y_t gramian} are well-defined;
	 	\item proving reasonable bounds on the moments of $Y_t$;
	 	\item
		assuring well-posedness as proposed in \cite{kazashi2025dynamical} for Lipschitz drift and  linear-growth bound of them (see e.g. \cite{kazashi2025dynamical});
	 	\item proving that the Gramian of $Y_t$ is always strictly positive in a matrix sense, and hence \eqref{eq: Y_t gramian} (which is the same equation of \eqref{eq:DLR-eq-Y}) is always well defined and the solution is global if it exists. For instance, this is the case if the diffusion is uniformly elliptic.
	 \end{itemize}
	 
	 \begin{Remark}[Approximation assumption]
		The strategy that we have followed in this section is based on a Euler-Maruyama discretization of the full-order solution \eqref{eq: full-order FE}. Notice that we have equivalent developments as long as another numerical method of weak order $1$ is employed to approximate the true solution in Section \ref{sec: full-order}.
	 \end{Remark}
	 
	 \begin{Remark}[Comparison with the equation in \cite{cao2018stochastic}]
	 	Equations \eqref{eq: U_t gramian}-\eqref{eq: Y_t gramian} are the same relations derived in \cite{cao2018stochastic}. Notice that this fact should not be surprising. Indeed, in \cite{cao2018stochastic} DO equations are derived as minimization of a measure constrained in a low-rank manifold. On the contrary, in this section we minimized the Gramian of a discretized surrogate, and hence, closely to minimize mean and covariance. If we conditioned all the aforementioned means and covariances over the the point $X_n^{\mathrm{DLRA}}$, then $X_{n+1}^{\mathrm{DLRA}}$ follows a Gaussian measure and, hence, it is completely defined by its first two moments.
	 \end{Remark}
	
	\subsection{DLRA obtained as minimization of mean and covariance}\label{sec: mean and cov DLRA}
	At the place of directly minimizing the Gramian, one can derive a DLRA solution via minimizing simultaneously the difference between the means and covariances of $X_{n+1}^{\mathrm{DLRA}}$ and of the element $X_{n+1}$ obtained by Euler-Maruyama approximation of the full order solution starting in the point $X_n^{\mathrm{DLRA}}$. As the process of derivation is similar to the discussion proposed for the Gramian in Section \ref{sec: min DLRA}, we sketch the computations of this approach here. Beyond relations \eqref{eq: matrix der}, we will use the following properties, too \cite{petersen2008matrix}: for all $a \in \mathbb{R}^p$, $b\in \mathbb{R}^m$ vectors and for $X \in \mathbb{R}^{p \times m}$ matrix one has
	\begin{equation*}
		\begin{aligned}
			&\frac{\partial X^{\top} a}{\partial X}
			= \frac{\partial a^{\top} X}{\partial X}
			= a, \qquad \frac{\partial a^{\top} X b}{\partial X}= ab^{\top},\\
			&\frac{\partial a^{\top} X^{\top} b}{\partial X}
			= ba^{\top}, \qquad
			\frac{\partial a^{\top} X a}{\partial X}= \frac{\partial a^{\top} X^{\top} a}{\partial X}
			= aa^{\top}.
		\end{aligned}
	\end{equation*}
	
	The mean of the surrogate $X_{n+1}$ defined in \eqref{eq: full-order FE} is 
		\begin{equation*}
		m_{n+1}:=\mathbb{E}[X_{n+1}] = \mathbb{E}[X_n^{\mathrm{DLRA}}] + \mathbb{E}[a(t, X_n^{\mathrm{DLRA}})] \Delta t_n,
	\end{equation*}
	and, hence, the centered $X_{n+1}$, i.e.\ $X_{n+1}$ minus its mean, namely $\mathring{X}_{n+1}$, is	
		\begin{equation*}
			\begin{aligned}
		\mathring{X}_{n+1} &= \mathring{X}_n^{\mathrm{DLRA}} + \left( a(t, X_n^{\mathrm{DLRA}}) - \mathbb{E}[a(t, X_n^{\mathrm{DLRA}})] \right) \Delta t_n + b(t, X_n^{\mathrm{DLRA}}) \Delta W_n \\
		& = \mathring{X}_n^{\mathrm{DLRA}} +  \mathring{a}_n \Delta t_n + b_n \Delta W_n,\\
		\end{aligned}
	\end{equation*}
	where the superscript $\mathring{f}$ indicates a random variable f minus its expectation, i.e.\ $\mathring{f}:= f- \mathbb{E}[f]$.
	
	The covariance $\mathring{C}_{n+1}$ at time $t_{n+1}$ is
	\begin{equation*}
		\begin{aligned}
		\mathring{C}_{n+1} &=\mathbb{E}[\mathring{X}_{n+1} \mathring{X}_{n+1}^\top] \\
		&= 
		\mathbb{E}[\mathring{X}_n^{\mathrm{DLRA}} (\mathring{X}_n^{\mathrm{DLRA}})^\top] + 
		\mathbb{E}[\mathring{X}_n^{\mathrm{DLRA}} \mathring{a}_n^\top] \Delta t_n + 
		\mathbb{E}[\mathring{a}_n (\mathring{X}_n^{\mathrm{DLRA}})^\top] \Delta t_n +
		\mathbb{E}[\mathring{a}_n \mathring{a}_n^\top] \Delta t_n^2 +
		\mathbb{E}[b_n b_n^\top] \Delta t_n.
		\end{aligned}
	\end{equation*}
	
	We compute the mean and covariance of $X_{n+1}^{\mathrm{DLRA}} = U_{n+1}^\top Y_{n+1}$, where $U_n \in \mathbb{R}^{k \times d}$ is a deterministic matrix with $U_n U_n^\top = I_{k \times k}$, i.e.\ $U_n$ has orthonormal rows, and $U_n \Delta U_n^\top = 0$, whereas $Y_n \in L^2(\Omega,\mathbb{R}^k)$ and $\Delta Y_n$ is adapted. The mean reads as
	\begin{equation*}
	m_{n+1}^{\mathrm{DLRA}} =	\mathbb{E}[X_{n+1}^{\mathrm{DLRA}}] = (U_n + \Delta U_n)^\top \left( \mathbb{E}[Y_n] + \mathbb{E}[\Delta Y_n] \right)
	\end{equation*}
	and the centered DLRA solution $\mathring{X}_{n+1}^{\mathrm{DLRA}}$ is defined as
	\begin{equation*}
		(\mathring{X}_{n+1}^{\mathrm{DLRA}}) := X_{n+1}^{\mathrm{DLRA}} - \mathbb{E}[X_{n+1}^{\mathrm{DLRA}}].
	\end{equation*}
	Furthermore, the covariance reads as
	\begin{equation*}
		\begin{aligned}
			\mathbb{E}[(\mathring{X}_{n+1}^{\mathrm{DLRA}}) (\mathring{X}_{n+1}^{\mathrm{DLRA}})^{\top}] =& (U_n + \Delta U_n)^\top \mathbb{E}[	\mathring{Y}_{n+1}  	\mathring{Y}_{n+1} ^\top] (U_n + \Delta U_n) \\
			\approx & U_n^\top \mathbb{E}[\mathring{Y}_n \mathring{Y}_n^\top] U_n + \Delta U_n^\top \mathbb{E}[\mathring{Y}_n \mathring{Y}_n^\top] U_n \\
			&+ U_n^\top \Delta C_{\mathring{Y}_n} U_n + U_n^\top  \mathbb{E}[\mathring{Y}_n \mathring{Y}_n^\top]  \Delta U_n \\
			\approx& U_n^\top  C_{\mathring{Y}_n} U_n + \Delta U_n^\top  C_{\mathring{Y}_n} U_n + U_n^\top \Delta C_{\mathring{Y}_n} U_n + U_n^\top   C_{\mathring{Y}_n}  \Delta U_n,
		\end{aligned}
	\end{equation*}
	where in the second line we employ a first order approximation as proposed in \eqref{eq: C DLRA}.
   
   In order to deal with the computation related to the minimization of $\| C_{n+1}^{\mathrm{DLRA}} - C_{n+1}\|_{\mathrm{F}}$ is useful to compute the following quantity
	\begin{equation*}
		\begin{aligned}
			\left(m_{n+1}^{\mathrm{DLRA}}\right)^\top m_{n+1}^{\mathrm{DLRA}} =& 
			(\mathbb{E}[Y_n] + \mathbb{E}[\Delta Y_n])^\top (U_n + \Delta U_n) (U_n^\top + \Delta U_n^\top)(\mathbb{E}[Y_n] + \mathbb{E}[\Delta Y_n]) \\
			=& (\mathbb{E}[Y_n] + \mathbb{E}[\Delta Y_n])^\top (I_{k \times k} + \Delta U_n \Delta U_n^\top)(\mathbb{E}[Y_n] + \mathbb{E}[\Delta Y_n]) \\
			=& \mathbb{E}[Y_n]^\top \mathbb{E}[Y_n] + \mathbb{E}[\Delta Y_n]^\top \mathbb{E}[Y_n] + \mathbb{E}[Y_n]^\top \mathbb{E}[\Delta Y_n] + \mathbb{E}[\Delta Y_n]^\top \mathbb{E}[\Delta Y_n] \\
			& + \mathbb{E}[Y_n]^\top \Delta U_n \Delta U_n^\top \mathbb{E}[Y_n] + \mathbb{E}[\Delta Y_n]^\top \Delta U_n \Delta U_n^\top \mathbb{E}[Y_n] \\
			& + \mathbb{E}[Y_n]^\top \Delta U_n \Delta U_n^\top \mathbb{E}[\Delta Y_n]  + \mathbb{E}[\Delta Y_n]^\top \Delta U_n \Delta U_n^\top \mathbb{E}[\Delta Y_n],
		\end{aligned}
	\end{equation*}
	where in the second line we employed the discrete gauge condition and the orthogonality of $U_n$. The derivative of the mean of the DLRA $m_{n+1}^{\mathrm{DLRA}}$ with respect to $\partial \mathbb{E}[\Delta Y_n]$ reads as
	\begin{equation*}
		\begin{aligned}
				\frac{\partial (m_{n+1}^{\mathrm{DLRA}})^\top m_{n+1}^{\mathrm{DLRA}} }{\partial \mathbb{E}[\Delta Y_n]} =\frac{\partial\mathbb{E}[X_{n+1}^{\mathrm{DLRA}}]^\top \mathbb{E}[X_{n+1}^{\mathrm{DLRA}}]}{\partial \mathbb{E}[\Delta Y_n]} =& 
			\mathbb{E}[Y_n] + \mathbb{E}[Y_n] + 2 \mathbb{E}[\Delta Y_n] \\
			&+ \Delta U_n \Delta U_n^\top \mathbb{E}[Y_n] + \Delta U_n \Delta U_n^\top\mathbb{E}[Y_n] + 2 \Delta U_n \Delta U_n^\top \mathbb{E}[\Delta Y_n] \\
			= & 2 \mathbb{E}[Y_n] + 2 \mathbb{E}[\Delta Y_n] + 2\Delta U_n \Delta U_n^\top \mathbb{E}[Y_n] + 2 \Delta U_n \Delta U_n^\top \mathbb{E}[\Delta Y_n]
		\end{aligned}
	\end{equation*}
	We compute now the derivative of the cross term
	\begin{equation*}
		\begin{aligned}
		m^\top_{n+1} m^{\mathrm{DLRA}}_{n+1} &= m^\top_{n+1} \left( U_n^\top + \Delta U_n^\top \right) \left( \mathbb{E}[Y_n] + \mathbb{E}[\Delta Y_n] \right),\\
		& = (U_n^\top \mathbb{E}[Y_n] +  \mathbb{E}[a_n]  \, \Delta t_n)^{\top} \left( U_n^\top + \Delta U_n^\top \right) \left( \mathbb{E}[Y_n] + \mathbb{E}[\Delta Y_n] \right) \\
		& = \mathbb{E}[Y_n]^{\top}  \mathbb{E}[Y_n] + \mathbb{E}[Y_n]^{\top} \mathbb{E}[\Delta Y_n] + \mathbb{E}[a_n^{\top}]  U_{n+1}^{\top} \mathbb{E}[Y_n]  \Delta t_n+  \mathbb{E}[a_n^{\top}] U_{n+1}^{\top}\mathbb{E}[\Delta Y_n] \Delta t_n
		\end{aligned}
	\end{equation*}
	with respect to the stochastic increment, namely
	\begin{equation*}
		\begin{aligned}
		\frac{\partial \, m^\top_{n+1} m^{\mathrm{DLRA}}_{n+1}}{\partial \, \mathbb{E}[\Delta Y_n]} &= \left( m^\top_{n+1} \left( U_n^\top + \Delta U_n^\top \right) \right)^\top = \mathbb{E}[Y_n] + U_n \mathbb{E}[a_n] \Delta t_n+ \Delta U_n \mathbb{E}[a_n] \Delta t_n,
		\end{aligned}
	\end{equation*}
thanks to the orthogonality of the rows of $U_n$ and the discrete gauge condition.
	Therefore, minimizing the mean in the Euclidean norm with respect to the stochastic increment $\Delta Y_n$ reads as
	\begin{equation}\label{eq: der mean DLRA}
		\begin{aligned}
			\frac{\partial \left| m_{n+1} - m^{\mathrm{DLRA}}_{n+1} \right|^2}{\partial \, \mathbb{E}[\Delta Y_n]} 
			=& \frac{ \left( -2 m^\top_{n+1} m^{\mathrm{DLRA}}_{n+1} + (m^{\mathrm{DLRA}}_{n+1})^{\top} m^{\mathrm{DLRA}}_{n+1} \right)}{\partial \, \mathbb{E}[\Delta Y_n]} \\
			=& 2 \, \mathbb{E}[Y_n] + 2 \, \mathbb{E}[\Delta Y_n] + 2 \, \Delta U_n \Delta U_n^\top \, \mathbb{E}[Y_n] + 2 \, \Delta U_n \Delta U_n^\top \, \mathbb{E}[\Delta Y_n] \\
			&- 2 \left( \mathbb{E}[Y_n] + U_n \mathbb{E}[a_n] \Delta t_n+ \Delta U_n \mathbb{E}[a_n] \Delta t_n \right)^\top\\
			=& 2 \, \mathbb{E}[\Delta Y_n] + 2 \, \Delta U_n \Delta U_n^\top \, \mathbb{E}[Y_n] - 2 \, \Delta U_n \Delta U_n^\top \, \mathbb{E}[\Delta Y_n] \\
			&\quad - 2 \, U_n \, \mathbb{E}[a_n] \, \Delta t_n - 2 \, \Delta U_n \, \mathbb{E}[a_n] \, \Delta t_n.
		\end{aligned}
	\end{equation}
	To find the minimum, we want that \eqref{eq: der mean DLRA} is null in a variational formulation, namely
	\begin{equation*}
		\left\langle \frac{\partial \left| m_{n+1} - m^{\mathrm{DLRA}}_{n+1} \right|^2}{\partial \, \mathbb{E}[\Delta Y_n]}, \, V \right\rangle = 0, \quad \forall V \in \mathbb{R}^{d \times k}, \ V \in \operatorname{Im}({U_n^{\top}}^\perp).
	\end{equation*}
	As the increment $\Delta U_n$ is at least of order $O(\Delta t_n)$, the in a first order of optimality, we have
	\begin{equation}\label{eq: Delta Y_n mean}
		\mathbb{E}[\Delta Y_n] = U_n \, \mathbb{E}[a_n] \, \Delta t_n.
	\end{equation}
	
	Now we minimize the Frobenius norm of the difference between the covariance of the full order solution and the DLRA one, similarly to the treatment proposed for the minimization of the Gramian approach. The covariance of the DLRA can be approximated in the following way in a first order of approximation in $\Delta t_n$
	\begin{equation*}
		\begin{aligned}
		\mathring{C}^{\mathrm{DLRA}}_{n+1} :=	\mathbb{E} \left[ \mathring{X}^{\mathrm{DLRA}}_{n+1} (\mathring{X}^{\mathrm{DLRA}}_{n+1})^\top \right] 
			=& (U_n + \Delta U_n)^\top \, \mathbb{E} \left[ \mathring{Y}_{n+1} \mathring{Y}_{n+1}^\top \right] \, (U_n + \Delta U_n) \\
			\approx & U_n^\top C_{\mathring{Y}_n} U_n + \Delta U_n^\top C_{\mathring{Y}_n} U_n + U_n^\top \Delta C_{\mathring{Y}_n} U_n + U_n^\top C_{\mathring{Y}_n} \Delta U_n \\
			=& U_n^\top \mathbb{E}[\mathring{Y}_{n}\mathring{Y}_{n}^\top] U_n + U_n^\top C_{\mathring{Y}_n} \Delta U_n + U_n^\top \Delta C_{\mathring{Y}_n} U_n +  \Delta U_n^\top \mathbb{E}[\mathring{Y}_{n}\mathring{Y}_{n}^\top] U_n.
		\end{aligned}
	\end{equation*}
	The covariance of the true solution is
	\begin{equation*}
		\begin{aligned}
			\mathring{C}_{n+1} = \mathbb{E}[\mathring{X}_{n+1} \mathring{X}_{n+1}^\top] 
			&= \mathbb{E}[\mathring{X}_n^{\mathrm{DLRA}} (\mathring{X}_n^{\mathrm{DLRA}})^{\top}] + \mathbb{E}[\mathring{X}_n^{\mathrm{DLRA}} \mathring{a}_n^\top] \Delta t_n \\
			&\quad + \mathbb{E}[\mathring{a}_{n} (\mathring{X}_n^{\mathrm{DLRA}})^{\top}] \Delta t_n + \mathbb{E}[\mathring{a}_{n} \mathring{a}_n^\top] \Delta t_n^2 + \mathbb{E}[b_n b_n] \Delta t_n.
		\end{aligned}
	\end{equation*}
	Putting these computations together, we found the following derivatives with respect to the  increment of the deterministic basis 
	\begin{equation*}
		\begin{aligned}
			\frac{\partial \, \mathrm{Tr} \left( \mathring{C}_{n+1}^{\mathrm{DLRA}} \left(\mathring{C}_{n+1}^{\mathrm{DLRA}}\right)^\top  \right)}{\partial \Delta U_n}
			&= 2 \, C_{\mathring{Y}_n} C_{\mathring{Y}_n} U_n + 4 \, C_{\mathring{Y}_n} C_{\mathring{Y}_n} \Delta U_n + 2 \, C_{\mathring{Y}_n} \Delta C_{\mathring{Y}_n} U_n \\
		\frac{\partial \mathrm{Tr} \left( \mathring{C}_{n+1}^{\mathrm{DLRA}} \mathring{C}_{n+1}^\top \right)}{\partial \Delta U_n} 
		&= 2 \, C_{\mathring{Y}_n} \, U_n \, \mathring{C}_{n+1}.
			\end{aligned}
	\end{equation*}
    Similarly for $\Delta Y_n$, we have
	\begin{equation*}
		\begin{aligned}
		\frac{\partial \text{Tr}\left( \mathring{C}_{n+1}^{\mathrm{DLRA}} (\mathring{C}_{n+1}^{\mathrm{DLRA}})^{\top} \right)}{\partial \Delta Y_n} 
		= 2 \, \Delta C_{\mathring{Y}_n} + 2 C_{\mathring{Y}_n}, \quad 
		\frac{\partial \text{Tr} \left( \mathring{C}_{n+1} (\mathring{C}_{n+1}^{\mathrm{DLRA}})^{\top} \right)}{\partial \Delta Y_n} 
		= U_n \, \mathring{C}_{n+1} \, U_n^{\top}.
		\end{aligned}
\end{equation*}
 Using the variational formulation as done with the minimization of the Gramian, one finds:
	\begin{equation*}
		\begin{aligned}
				\Delta C_{\mathring{Y}_n} = U_n \, \mathring{C}_{n+1} \, U_n^{\top} - C_{\mathring{Y}_n},
				\quad \Delta U_n = C_{\mathring{Y}_n}^{-1} \, U_n \, \mathring{C}_{n+1}\, P_{U_n}^\perp.
			\end{aligned}
	\end{equation*}
We use the constaints on $\Delta C_{\mathring{Y}_n}$ to obtain an equation for $\Delta \mathring{Y}_n$. From one side we have,
	\begin{equation}\label{eq: Delta C_y 1}
		\begin{aligned}
			\Delta C_{\mathring{Y}_n} =& U_n \Big( 
			\mathbb{E}[\mathring{X}_n^{\mathrm{DLRA}} (\mathring{X}_n^{\mathrm{DLRA}})^{\top}] + \mathbb{E}[\mathring{X}_n^{\mathrm{DLRA}} (\mathring{a}_{n})^{\top}] \Delta t_n 
			+ \mathbb{E}[\mathring{a}_{n} (\mathring{X}_n^{\mathrm{DLRA}})^{\top}] \Delta t_n \\
			& + \mathbb{E}[\mathring{a}_{n} (\mathring{a}_{n})^{\top}] \Delta t_n^2 
			+ \mathbb{E}[b_n b_n^{\top}] \Delta t_n 
			\Big) U_n^{\top} - C_{\mathring{Y}_n} \\
			=& \mathbb{E}[\mathring{Y}_n \mathring{Y}_n^{\top}] \Delta t_n 
			+ U_n \mathbb{E}[\mathring{a}_{n} \mathring{Y}_n^{\top}] \Delta t_n 
			+ U_n \mathbb{E}[\mathring{a}_{n} (\mathring{a}_{n})^{\top}] U_n^{\top} \Delta t_n^2 
			+ U_n \mathbb{E}[b_n b_n^{\top}] U_n^{\top} \Delta t_n,
		\end{aligned}
	\end{equation}
	whereas, on the other hand, one has
	\begin{equation}\label{eq: Delta C_y 2}
		\Delta C_{\mathring{Y}_n} = \mathbb{E}[\mathring{Y}_n \mathring{Y}_n^{\top}] + \mathbb{E}[\Delta \mathring{Y}_n \mathring{Y}_n^{\top}] + \mathbb{E}[\Delta \mathring{Y}_n \Delta \mathring{Y}_n^{\top}].
	\end{equation}
	Matching \eqref{eq: Delta C_y 1} with \eqref{eq: Delta C_y 2}, one obtains $\Delta Y_n = U_n \, \mathring{a}_{n} \, \Delta t_n + U_n \, b_n \, \Delta W_n$
	and, hence,
	\begin{equation*}
		\mathring{Y}_{n+1} = \mathring{Y}_n + \Delta \mathring{Y}_n = \mathring{Y}_n + U_n \, \mathring{a}_{n} \, \Delta t_n + U_n \, b_n \, \Delta W_n.
	\end{equation*}
	As the average of $Y_{n+1}$ reads as $\mathbb{E}[Y_{n+1}] = \mathbb{E}[Y_n] + U_n \, \mathbb{E}[a_{n}] \Delta t_n$,
	then one retrieves the usual stochastic update for $Y_{n+1}$
		\begin{equation}\label{eq: Y_n mean cov}
		Y_{n+1} = Y_n + U_n \, a_{n} \Delta t_n + U_n \, b_n \, \Delta W_n.
	\end{equation}
	For the increment $\Delta U_n$ of the deterministic basis, one gets
	\begin{equation*}
		\begin{aligned}
			\Delta U_n &= C_{\mathring{Y}_n}^{-1} \mathbb{E}[\mathring{Y}_n (\mathring{a}_{n})^{\top}] P_{U_n}^\perp \Delta t_n 
			+ C_{\mathring{Y}_n}^{-1} U_n \mathbb{E}[\mathring{a}_{n} (\mathring{a}_{n})^{\top}] P_{U_n}^\perp \Delta t_n^2 + C_{\mathring{Y}_n}^{-1} U_n \mathbb{E}[b_n b_n^{\top}] P_{U_n}^\perp \Delta t_n \\
			&\approx C_{\mathring{Y}_n}^{-1} \mathbb{E}[\mathring{Y}_n (\mathring{a}_{n})^{\top}] P_{U_n}^\perp \Delta t_n 
			+ C_{\mathring{Y}_n}^{-1} U_n \mathbb{E}[b_n b_n^{\top}] P_{U_n}^\perp \Delta t_n,
		\end{aligned}
	\end{equation*}
	where in the last line we discarded the terms with order higher than $\Delta t_n$, as done in the case of the Gramian minimization for \eqref{eq: U_n gramian}. Therefore, the equation for $U_{n}$ is 
	\begin{equation}\label{eq: U_n mean cov}
		\begin{aligned}
			U_{n+1} &= U_n + \Delta U_n \\
			&= U_n + C_{\mathring{Y}_n}^{-1} \mathbb{E}[\mathring{Y}_n (\mathring{a}_{n})^{\top}] P_{U_n}^\perp \Delta t_n 
			+ C_{\mathring{Y}_n}^{-1} U_n \mathbb{E}[b_n b_n^{\top}] P_{U_n}^\perp \Delta t_n.
		\end{aligned}
	\end{equation}
	Notice that, \eqref{eq: U_n mean cov} differs from \eqref{eq: U_n gramian} also by the presence of centered terms $C_{\mathring{Y}_n}, \mathring{Y}_n,$ and $\mathring{a}_{n}$.
	
	As done in Section \ref{sec: continuous time limit}, we suppose that $(U_n)_n$ and $(Y_n)_n$ are sequences approximating a deterministic function $(U_t)_{t\geq 0}$ and a stochastic one $(Y_t)_{t\geq 0}$, respectively, in the time mesh $(t_n)_n$, where $U$,$Y$ possess properties $1.,2.,3.$ of Definition \ref{def: D0 sol}. We define the equations that $U$ and $Y$ satisfy as the limit 
	for $\Delta t_n \to 0$ of \eqref{eq: U_n mean cov} and \eqref{eq: Y_n mean cov}, namely
	\begin{align}
			\mathrm{d}U_t &= \mathring{C}_{Y_t}^{-1} \mathbb{E}[\mathring{Y}_t \mathring{a}(t,X_t)^\top] P_{U_t}^\perp \mathrm{d}t+ \mathring{C}_{Y_t}^{-1} U_t \mathbb{E}[b(t,X_t) b(t,X_t)^\top] P_{U_t}^\perp \mathrm{d}t \label{eq: U_t mean cov},\\
		\mathrm{d}Y_t &= U_t a(t,X_t)\,\mathrm{d}t + U_t b(t,X_t)\,\mathrm{d}W_t \label{eq: Y_t mean cov}
	\end{align}
	where $\mathring{Y}_t$ is the centered stochastic basis, i.e.\ $\mathring{Y}_t:= Y_t - \mathbb{E}[Y_t]$, $\mathring{C}_{Y_t}:=\mathbb{E}[\left(Y_t - \mathbb{E}[Y_t]\right)\left(Y_t - \mathbb{E}[Y_t]\right)^{\top}]$ is the covariance of $Y_t$, and $\mathring{a}(t,X_t):= a(t,X_t) - \mathbb{E}[a(t,X_t)]$ is the centered drift.

    The existence and uniqueness of solutions for \eqref{eq: U_t mean cov} and \eqref{eq: Y_t mean cov} follows completely similarly to the treatment discussed for \eqref{eq: U_t gramian} and \eqref{eq: Y_t gramian}.
	
	\begin{Remark}[Consistency of DLRA approaches]\label{rmk: consistency}
		As in this section, we derived different DO equations to the ones obtained in \cite{kazashi2025dynamical}, it is reasonable to query which surrogate dynamics is the best. Unfortunately, the consistency argument developed in \cite{kazashi2025dynamical} does not give an answer. Indeed, if 
		$U_{t}^{\top}\beta_{t} =b(t,U_{t}^{\top}Y_{t})$
		for some $\beta_{t}\in\mathbb{R}^{R\times m}$ progressively measurable
		and with continuous paths almost surely, then
		\begin{equation*}
			C_{Y_t}^{-1} U_t \mathbb{E}[b_t b_t^{\top}] P_{U_t}^\perp  = C_{Y_t}^{-1} \mathbb{E}[\beta_{t}\beta_{t}^{\top} U_{t}] P_{U_t}^\perp = 0,
		\end{equation*}
		which means that the additional term is null if the surrogate to approximate $X_t^{\mathrm{true}}$ and, hence, the diffusion term $b$, is always low-rank with image in $\mathrm{Ran}(U_t)$. Therefore, \eqref{eq: U_t gramian} and \eqref{eq:DLR-eq-U} differ only when the solution to approximate is characterized by drift and diffusion that span in different directions.
		The same discussion holds for the DO equations in \cite{cao2018stochastic}.
    \end{Remark}
    \begin{Remark}[Three-term DLRA]
    	One can follow a similar treatment of this section to derive a three-term DLRA composed by the triplet $(m_n^{\mathrm{DLRA}}, U_n, Y_n)$. Notice that in that case the mean would follow the same equation as the mean of the true solution.
    \end{Remark}
	\subsection{A Possible Stable Numerical Algorithm}
	By construction, equations \eqref{eq: U_n gramian} and \eqref{eq: Y_n gramian} represent a possible discretization of \eqref{eq: U_t gramian} and \eqref{eq: Y_t gramian}, respectively.
	 The relation \eqref{eq: Y_n gramian} is the same equation of a DLR Projector Splitting for SDE \cite[Algorithm 3]{kazashi2026dynamicalpartI} and, hence, one can expect the boundedness of its second moments. We can define a staggered method that approximates  \eqref{eq: Y_t gramian} and \eqref{eq: U_t gramian} in the same fashion of \cite[Algorithm 3]{kazashi2026dynamicalpartI}, which is numerically stable and convergent. Namely, we have:
	 \begin{enumerate}
	 	\item Compute  $\tilde{Y}_{n+1}$ and $\tilde{U}_{n+1}$ as solutions of ($a_n=a(t_n,U_{n}^{\top}Y_n),$ $b_n=b(t_n,U_{n}^{\top}Y_n)$)
	 	\begin{equation*}
	 		\begin{aligned}
	 			 \tilde{Y}_{n+1} &= Y_n + U_n a_n \Delta t_n +  U_n b_n \Delta W_n\\
	 			 C_{\tilde{Y}_{n+1}} \tilde{U}_{n+1} &=  C_{\tilde{Y}_{n+1}} U_n + \left(\mathbb{E}\left[\tilde{Y}_{n+1} a_n^{\top} \right] + U_n \mathbb{E}[b_n b_n^{\top}]\right) \left(I_{d \times d} - P_{U_n} \right) \Delta t_n
	 		\end{aligned}
	 	\end{equation*}
	 	\item Reorthonormalize the deterministic modes: find  $({U}_{n+1}, {Y}_{n+1})$ such that:
	 	\begin{equation*}
	 	U_{n+1}^{\top} Y_{n+1} = \tilde{U}_{n+1}^{\top} \tilde{Y}_{n+1}, \quad {U}_{n+1}{U}_{n+1}^{\top} = I_{d\times d}.
	 	\end{equation*}
	 	with $(U_{n+1}^{\top}, R) =\texttt{QR}(\tilde{U}_{n+1}^{\top})$ and $Y_{n+1} = R \tilde{Y}_{n+1}$, where $\texttt{QR}$ is the QR decomposition: for all rectangular matrices $A \in \mathbb{R}^{n \times k}$ one has $(Q,R) =\texttt{QR}(A)$ with $Q \in \mathbb{R}^{n \times k}$ with orthogonal columns and $R \in \mathbb{R}^{k \times k}$. 
	 	\item Set $X_{n+1}^{\mathrm{DLRA}}:=U_{n+1}^{\top} Y_{n+1}$.
	 \end{enumerate}
	 Then, the surrogate update of $X_{n+1}^{\mathrm{DLRA}}$ reads as follows
	\begin{equation}\label{eq: X_n eq}
	\begin{aligned}
			X_{n+1}^{\mathrm{DLRA}}  =& X_n^{\mathrm{DLRA}} + P_{U_n} a_n \Delta t_n + P_{U_n} b_n \Delta W_n \\
			&+ P_{U_n}^\perp P_{\widetilde{Y}_{n+1}} [a_n] \Delta t_n + P_{U_n}^\perp \mathbb{E}[b_n b_n^\top]  U_n^{\top} C_{ \widetilde{Y}_{n+1}}^{-1} \widetilde{Y}_{n+1} \Delta t_n.
		\end{aligned}
	\end{equation}

     Even though \eqref{eq: X_n eq} is characterized by the inverse of the Gramian, the boundedness of $\mathbb{E}[|X_{n+1}^{\mathrm{DLRA}}|^2]$ of the \cite[Algorithm 3]{kazashi2026dynamicalpartI} still holds with this DO formulation.
   Indeed, we have the following chain of equivalences
   \begin{equation}\label{eq: chain of eq P}
   	\begin{aligned}
   		&P_{U_n}^\perp \mathbb{E}[b_n b_n^\top] U_n^{\top} C_{ \widetilde{Y}_{n+1}}^{-1} \widetilde{Y}_{n+1}\Delta t_n \\
   		= & P_{U_n}^\perp \mathbb{E}[b_n  \Delta t_n  b_n^\top U_n^{\top} ] C_{ \widetilde{Y}_{n+1}}^{-1} \widetilde{Y}_{n+1} \\
   		= & P_{U_n}^\perp \mathbb{E}[b_n  \Delta W_n \Delta W_n^{\top}  b_n^\top U_n^{\top} ] C_{ \widetilde{Y}_{n+1}}^{-1} \widetilde{Y}_{n+1} \\ 
   		= & P_{U_n}^\perp \mathbb{E}[b_n  \Delta W_n (Y_n+ U_n a_n \Delta t_n +U_n b_n \Delta W_n - Y_n - U_n a_n \Delta t_n)^{\top} ] C_{ \widetilde{Y}_{n+1}}^{-1} \widetilde{Y}_{n+1} \\ 
   		= & P_{U_n}^\perp \mathbb{E}[b_n  \Delta W_n (Y_n+ U_n a_n \Delta t_n +U_n b_n \Delta W_n)^{\top} ] C_{ \widetilde{Y}_{n+1}}^{-1} \widetilde{Y}_{n+1}  =   P_{U_n}^\perp P_{\widetilde{Y}_{n+1}} [b_n \Delta W_n ],\\
   	\end{aligned}
   \end{equation}
   where to pass from the second-to-last to the last line we use the independence of the increments $\Delta W_n$. Therefore, similarly to the result proved in \cite[Lemma 6.1]{kazashi2026dynamicalpartI} this algorithm gives a discretized solution with bounded second moment, i.e.\ $\sup_{1\leq n\leq N} \mathbb{E}[|X_{n}^{\mathrm{DLRA}}|^2] < \infty$, under linear-growth bound assumption on drift $a$ and diffusion $b$.
   
  One can prove convergence of this algorithm with respect to the true solution up to a standard approximation assumption. Details of the needed procedure are very similar to the one developed in \cite[Sections 5 and 6]{kazashi2026dynamicalpartI}.
   
   	\section{New Approaches: a Stratonovich derivation}\label{sec: strat dlra}
	In this section, we try to investigate another possible surrogate approximation of SDEs that resembles a DLRA formulation which exploits the Stratonovich formalism and the chain rule formula. Indeed, unlike Itô formalism, in this case the stochastic integral is built so that the chain rule, defined as in ordinary calculus, holds. This setting is beneficial to directly build stochastic processes constrained on manifolds, at the price of losing martingale properties for the stochastic integral.
	
	In \cite{koch2007dynamical}, Koch and Lubich derived DLRA for ODEs, namely for a master matrix differential equation $\dot{X}^{\mathrm{true}}_t = \mathrm{F}(X^{\mathrm{true}}_t)$ with $X^{\mathrm{true}} \in \mathbb{R}^{n \times m}$ for all $t\geq 0$ and $\mathrm{F}: \mathbb{R}^{n \times m} \to \mathbb{R}^{n \times m}$, as a minimization problem. Indeed, they were looking for a matrix $X_t \in \mathbb{R}^{n \times m}$ of rank $k$ satisfying
	\begin{equation}\label{eq: dlra}
		X_t = \min\limits_{Z_t \in \mathcal{M}_{k}} \| \dot{Z}_t - \mathrm{F}(Z_t)\|_{\mathrm{F}}, \text{ for all } t\geq 0,
	\end{equation}
	where $\mathcal{M}_{k}$ denotes the manifold of rank-$k$ $n \times m$-matrices and $\mathcal{T}_{\mathrm{X}_t}\mathcal{M}_{k}$ the tangent space of the manifold $\mathcal{M}_{k}$ at the point $\mathrm{X}_t$.
	Formally, problem \eqref{eq: dlra} describes the best low-rank surrogate whose difference between its derivative and the right-hand side of the studied equation was the minimum possible in a given norm (e.g. as the Frobenius one in \eqref{eq: dlra}). In the same article it was proven that this minimization problem has a direct geometrical equivalent interpretation. Indeed, \eqref{eq: dlra} can be rewritten in the following form
		\begin{equation}\label{eq: dlra geom}
		X_t  \in \mathcal{M}_{k} \text{ such that } X_t = P_{X_t} \mathrm{F}(X_t),
	\end{equation}
	where $P_{X_t}$ denotes the orthogonal projector onto $\mathcal{T}_{\mathrm{X}_t}\mathcal{M}_{k}$. Relation \eqref{eq: dlra geom} means that
	the derivative of the DLRA $X_t$ for ODE is equal to the orthogonal projection in the point $X_t$ onto the tangent space of a manifold of low-rank functions applied to the right hand side of the ODE, computed in the same point $X_t$. 
	
	To give a practical visualization, the DLRA $X_t$ shows the following behavior between times $t_1$ and $t_2$:
	 \begin{multicols}{2}
	 		\begin{center}
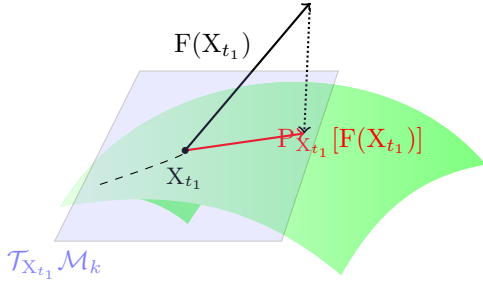

	 			\begin{tikzpicture}[scale=0.75]	
	 				\draw[fill=green!50, draw=none, shift={(0.1, 0.6)},scale=1.5]
	 				(0, 0) to[out=20, in=140] (1.5, -0.2) to [out=60, in=160]
	 				(5, 0.5) to[out=130, in=60]
	 				cycle
	 				node[anchor=west, green, right=2.7cm, yshift=2em]{$\mathcal{M}_k$};
	 				
	 				\shade[thin, left color=green!10, right color=green!50, draw=none,
	 				shift={(0.1, 0.6)},scale=1.5]
	 				(0, 0) to[out=10, in=140] (3.3, -0.8) to [out=60, in=190] (5, 0.5)
	 				to[out=130, in=60] cycle;
	 				
	 				\draw[->, red, thick, shift={(0.1, 0.6)}] (2.2, 1) -- ++(2.1, .3) node[midway, below=0.2cm,right=0.3cm]{$\mathrm{P}_{\mathrm{X}_{t_1}}[\mathrm{F}(\mathrm{X}_{t_1})]$};
	 				
	 				\filldraw[black, shift={(0.1, 0.6)}] (2.2,1) circle (1.5pt) node[anchor=north, below=0.1cm]{\small $\mathrm{X}_{t_1}$};
	 				
	 				\draw[->, thick, shift={(0.1, 0.6)}] (2.2, 1) -- ++(2.2, 2.6) node[midway, above, xshift=-1.3em, yshift=.3em]{$\mathrm{F}(\mathrm{X}_{t_1})$};
	 				
	 				\draw[->, thick, shift={(0.1, 0.6)}, densely dotted] (4.4, 3.6) -- (4.3,1.3);
	 				\filldraw[fill=blue!40!white, opacity=.2] (0,0) -- (4,0) -- (5,3.) -- (1.5,3.) -- cycle
	 				node[anchor=north, opacity=.5, blue]{$\mathcal{T}_{\mathrm{X}_{t_1}}\mathcal{M}_{k}$};
	 				\draw[shift={(0.1, 0.6)}, thin, dashed] (.7,.4) .. controls (1,.5) and (2.5,1) .. (2.2,1);
	 		\end{tikzpicture}
	 		\captionof{figure}{DLRA $X$ at time $t_1$}
	 			\label{fig: dlra t_1}
	 	\end{center}
	 		\columnbreak
	 			\begin{center}
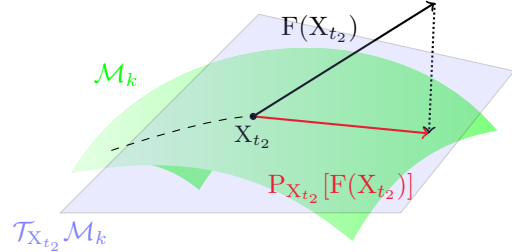

	 			\begin{tikzpicture}[scale=0.75]	
	 				\draw[fill=green!50, draw=none, shift={(0.2, 0.7)},scale=1.5]
	 				(0, 0) to[out=20, in=140] (1.5, -0.2) to [out=60, in=160]
	 				(5, 0.5) to[out=130, in=60]
	 				cycle
	 				node[anchor=west, green, yshift=2em]{$\mathcal{M}_k$};
	 				
	 				\shade[thin, left color=green!10, right color=green!50, draw=none,
	 				shift={(0.2, 0.7)},scale=1.5]
	 				(0, 0) to[out=10, in=140] (3.3, -0.8) to [out=60, in=190] (5, 0.5)
	 				to[out=130, in=60] cycle;
	 				
	 				\draw[->, red, thick, shift={(0.2, 0.7)}] (3.2, 1) -- ++(3.1, -.3) node[midway, below=0.5cm]{$\mathrm{P}_{\mathrm{X}_{t_2}}[\mathrm{F}(\mathrm{X}_{t_2})]$};
	 				
	 				\filldraw[black, shift={(0.2, 0.7)}] (3.2,1) circle (1.5pt) node[anchor=north]{\small $\mathrm{X}_{t_2}$};
	 				
	 				\draw[->, thick, shift={(0.2, 0.7)}] (3.2, 1) -- ++(3.2, 2) node[midway, above, xshift=-.9em, yshift=.3em]{$\mathrm{F}(\mathrm{X}_{t_2})$};
	 				
	 				\draw[->, thick, shift={(0.2, 0.7)}, densely dotted] (6.4, 3) -- (6.3,.7);
	 				\filldraw[fill=blue!40!white, opacity=.2] (0,0) -- (6,0) -- (8,2.5) -- (3.5,3.5) -- cycle
	 				node[anchor=north, opacity=.5, blue]{$\mathcal{T}_{\mathrm{X}_{t_2}}\mathcal{M}_{k}$};
	 				\draw[shift={(0.2, 0.7)}, thin, dashed] (.7,.4) .. controls (1,.5) and (2.5,1) .. (3.2,1);
	 			\end{tikzpicture}
	 			\captionof{figure}{DLRA $X$ at time $t_2$}
	 			\label{fig: dlra t_2}
	 				\end{center}´
	 	\end{multicols}
	On Figure \ref{fig: dlra t_1}, the given dynamics tends to exit the manifold $\mathcal{M}_k$ due to the force term $F(X)$ , but the solution is projected back to $\mathcal{M}_k$ thanks to the orthogonal projector $\mathrm{P}_{\mathrm{X}(t_1)}$. From this projection the point $X(t_2)$ is obtained after integrating and then the whole procedure is iterated as seen in Figure \ref{fig: dlra t_2}.
		
	The geometrical interpretation of DLRA requires the differentiability of the solution and, hence, a well-definition of the tangent space in a point of a given manifold. 
	Itô SDEs do not possess the derivability of the solution with respect to time. That is why in \cite{kazashi2025dynamical} a stochastic calculus strategy via Itô's formula is pursued to derive the DO equations. However, one can observe that considering the Stratonovich SDE interpretation would allow to use the chain rule on coefficients, as the standard differential calculus, allowing to have a similar setting to the one of \cite{koch2007dynamical}. This other formalism can be pursued at the price of losing martingale properties \cite{protter2012stochastic}. 
	
	The Stratonovich interpretation is at the base of the definition of stochastic calculus on manifold. Indeed, the common way to define surrogate processes is through a projection of the drift and the diffusion onto a finite dimensional manifold (see e.g.\ \cite{hsu2002stochastic,ito1950stochastic}). This is not the case for the DLRA expressed in \cite{kazashi2025dynamical}, where the projection onto the corange of the process is done in $L^2(\Omega)$. Therefore, it makes sense to use Stratonovich calculus to try to derive equations à la DLRA for a time-dependent surrogate.
	
	With these considerations in mind, we want to define a reasonable DLRA-type approximation using the standard stochastic calculus on manifold. We streamline our idea here before giving a more detailed discussion. 
	To apply the machinery previously described, we consider an ensemble of $M$ realizations of an SDE, which describes a noisy particle system. The evolution of the whole ensemble can be described by a matrix SDE living in some finite dimensional space, where each column represents a realization of the previous standard SDE. On this setting we can build our DLRA using the stochastic calculus on manifold, i.e.\ by constraining the latter matrix SDE to live on alow-dimensional manifold via projection of the drift and diffusion. Finally, if we sample independently these $M$ realizations, under some standard conditions, for a large number $M$ we expect that the empirical measure of these $M$ realizations converges to the measure of the mean-field SDE (associated to the noisy particle system). Therefore, we will derive standard-type DLRA equations for SDEs by considering the mean-field limit of the DLRA for the aforementioned matrix SDE, i.e.\ by (formally) seeing to which surrogate this matrix DLRA system converges for large $M$.
	
	\subsection{Construction of the DLRA particle system}
	We consider a discretization of a candidate rank-$k$ DLRA $X$ over the stochastic space using a Monte-Carlo method with $M$ samples. $X$ is therefore a matrix of dimension $d \times M$, where $d$ is the physical dimension, of rank $k$. We assume that each sample is described by a Stratonovich SDE where each drift and diffusion is projected in the common tangent space at the discretized point $X$ of the manifold of rank-$k$ matrix. Then, we translate this matrix SDE system in a Itô form and we consider the (formal) limit for the number of samples $M$ which goes to infinity. Assuming the existence of this limit, we consider this limit as our new DLRA formulation. We present our strategy more in details hereafter. To work into the Stratonovich framework, we assume that our diffusion $b$ is differentiable with respect to the spatial coordinate.
	
	For all $t\geq 0$, we consider a stochastic process $(X_t^M)_{t \geq 0}$, with $X_t^M \in \mathbb{R}^{d \times M}$ for each $t\geq0$, where the superscript $M$ can be thought as the number of
	samples in the Monte-Carlo discretization of a process in $L^2(\Omega, \mathbb{R}^d)$. With the notation $X^{(M,i)} \in \mathbb{R}^d$  we denote the $i$-th
	realization of $X_t^M$ , for $i=1,\ldots,M$, i.e.\ its $i$-th column.
	
	In this section, let us consider the following probability space $\left( \Omega, \mathcal{F}, \mathbb{P}, (\mathcal{F}_t)_{t \geq 0} \right)$, where $\Omega$, the probability domain, is a Polish space, $\mathcal{F}$ is a $\sigma$-algebra on $\Omega$, $\mathbb{P}$ is a measure of probability on $\Omega$ and $(\mathcal{F}_t)_{t \geq 0}$ is a standard filtration on the probability space $\left(\Omega, \mathcal{F}, \mathbb{P}\right)$. We suppose that our surrogate $X^{M,i}$ approximates the problem \eqref{eq:SDE-int} for large $M$.
	\begin{Remark}
		The choice of asking for $\Omega$ to be Polish is made to guarantee the convergence of the empirical measure given by the Monte-Carlo method applied to $X^{(M,i)}_{i=1,\dots,M}$ to the true measure in $L^2(\Omega)$ as a consequence of Varadarajan theorem \cite[Theorem 11.4.1]{dudley2018real}. A necessary condition under which the result holds is uniform boundedness of the second moment of $X^{(M,i)}$, i.e.\ $\sup_{i=1,\dots,M} \mathbb{E}[|X^{(M,i)}|^2] < \infty$. This last property can be proven under usual linear-growth bound for standard DLRA equations (see \cite{kazashi2025dynamical}).
	\end{Remark}
	\begin{Remark}
		Asking for the probability space $\Omega$ to be Polish is not such a restrictive assumption. 
			Indeed, such space can be constructed as follows. 
			First, consider the Wiener Space $(C_0\left([0,\infty);\mathbb{R}^m\right),\mu)$, where $C_0\left([0,\infty);\mathbb{R}^m\right)$ is the space of all continuous function $f: [0,\infty) \to \mathbb{R}^m$, with $f(0)=0$, endowed with the metric $\mathfrak{d}$ defined as
			\[\mathfrak{d}(f,g) = \sum_{n=1}^{\infty} \frac{1}{2^n}\sup\limits_{0 \leq t \leq n} |f(t) - g(t)| \wedge 1, \quad f,g \in C_0\left([0,\infty);\mathbb{R}^m\right),\]and $\mu$ is the Wiener measure. 
			The space $\left(C_0\left([0,\infty);\mathbb{R}^m\right), \mathfrak{d}\right)$ is complete and separable \cite{schilling2021brownian}.
			Moreover, assume that $X_{0}^{\mathrm{true}} \sim \nu$, where $\nu : \mathbb{R}^d \to [0,1]$ is a probability measure on $\left(\mathbb{R}^d, \mathcal{B}(\mathbb{R}^d)\right)$ where  $\mathcal{B}(\mathbb{R}^d)$ is the Borel $\sigma$-algebra of $\mathbb{R}^d$. 
			With the Euclidean norm, $\mathbb{R}^d$ is a Polish space. 
			Then, we can build our probability space $(\Omega,\mathcal{F},\mathbb{P})$ as follows. 
			As $\Omega$ we take the product metric space between the classical Wiener Space and $\mathbb{R}^d$, i.e.\ $\Omega :=C_0\left([0,\infty);\mathbb{R}^m\right) \times \mathbb{R}^d$. 
			By construction, the product space $\Omega$ is still Polish. 
			With respect to its product metric topology we take the Borel $\sigma$-algebra $\mathcal{F}_{\text{pre}}\subset 2^\Omega$, which is identical to the product $\sigma$-algebra of the two $\sigma$-algebras because of the separability of the spaces under discussion, and the product measure $\mathbb{P}_{\text{pre}}$ of $\mu$ and $\nu$ on $\mathcal{F}_{\text{pre}}$. Finally we take $(\Omega,\mathcal{F},\mathbb{P})$ as the completion of $(\Omega,\mathcal{F}_{\text{pre}},\mathbb{P}_{\text{pre}})$.
	\end{Remark}
	
	We suppose that for any $M$ our rank-$k$ surrogate $X^M$ is composed by the product of a deterministic and a stochastic basis $U^M$ and $Y^M$, respectively, namely
	$$X_t^M = (U^M_t)^{\top} Y_t^M \in \mathbb{R}^{d \times M},$$
    where for all $t\geq0$, $U^M_t \in \mathbb{R}^{k \times d}$ has orthogonal rows, i.e.\ $U^M_t(U^M_t)^\top = I_{k \times k}$, and $U^M_t $ satisfies a gauge condition, namely $\delta U^M_t (U^M_t)^\top = 0$, with $\delta U^M_t $ variation of $U^M_t $, and $Y^M_t \in \mathbb{R}^{k \times M}$ has full rank $k$. 

	For all $t$ and $M$,
	one can define the map $X_t^M : (\mathbb{R}^M, \langle \cdot, \cdot \rangle_M) \to (\mathbb{R}^d, \langle \cdot, \cdot \rangle)$, where  $\langle \cdot, \cdot \rangle$ is the Euclidean scalar product and 	\[
	\langle x,y \rangle_M = \frac{1}{M} \langle x,y \rangle
	\]
	is the scalar product derived from the empirical measure defined by
	the Monte--Carlo estimator, as
	\begin{equation*}
		\mathbb{R}^M \ni v \mapsto X_t^Mv \in \mathbb{R}^d.
	\end{equation*}
	By construction, the range of $X_t^M$ is defined in $(\mathbb{R}^d, \langle \cdot, \cdot \rangle)$
	and its corange in $(\mathbb{R}^d, \langle \cdot, \cdot \rangle_M)$.
	Let us define
	$\mathcal{M}_k(\mathbb{R}^{d \times M}) = \bigl\{ A \in \mathbb{R}^{d \times M} \,\big|\, \operatorname{rank}(A) = k \bigr\},$ i.e.\
	the manifold whose elements are matrices  in $\mathbb{R}^{d \times M}$ of rank $k$. Having these observations in mind, we now try to derive useful quantities for a general rank-$k$ element of the manifold.
	
	Given a generic element $X \in \mathcal{M}_k(\mathbb{R}^{d \times M})$,
	consider its decomposition $X = U^\top Y$, with $U\in \mathbb{R}^{k \times d}$ having orthonormal rows and satisfying a gauge condition and $Y\in \mathbb{R}^{k \times M}$ of full rank. Then, its variation can be written as
	$$\delta X = \delta U^\top Y + U^\top \delta Y,$$
	and the tangent space in the point $X$ of $ \mathcal{M}_k(\mathbb{R}^{d \times M})$ is defined as
	\begin{equation}\label{eq: Tan X}
		T_X \mathcal{M}_k(\mathbb{R}^{d \times M}) = 
	\Bigl\{ \delta X = \delta U^\top Y + U^\top \delta Y \,\big| \,
	U \delta U^\top = 0,\text{ with } U,\delta U \in \mathbb{R}^{k \times d},\;
	Y,\delta Y \in \mathbb{R}^{k \times M}, \ X = U^\top Y \Bigr\}.
	\end{equation}	
	Notice that $T_X \mathcal{M}_k(\mathbb{R}^{d \times M})$ can be interpreted as a subset of the ambient space $\mathbb{R}^{d \times M}$, and, hence, all the well-defined operations in $\mathbb{R}^{d \times M}$ are inherited here.
	
	\subsection{The orthogonal projector onto $T_X \mathcal{M}_k(\mathbb{R}^{d \times M})$}\label{sec: orth proj}
    In order to obtain differential equation for our DLRA, we want to obtain a explicit formula for the orthogonal projector
	$$P(X) : \mathcal{M}_k(\mathbb{R}^{d \times M}) \to T_X \mathcal{M}_k(\mathbb{R}^{d \times M}),
	\qquad
	X \mapsto P_X,$$
	i.e.\ $P(X) := P_X$ for the sake of notation.
	In the Stratonovich formulation, derivatives of the diffusion appear (providing that this diffusion is differentiable). In the DLRA setting, we then need an expression for the derivative of the projector $P_X$. In order to derive the expression of this differential, we need to analyze the decomposition of the tangent space at the point $X = U^\top Y$
	into its vertical and horizontal components, namely
	$$V_X \mathcal{M}_k(\mathbb{R}^{d \times M}) =
	\Bigl\{ U^\top \delta Y \,\big|\, \delta Y \in \mathbb{R}^{k \times M}, X = U^\top Y \Bigr\},$$
	and
	$$H_X \mathcal{M}_k(\mathbb{R}^{d \times M}) =
	\Bigl\{ \delta U^\top Y \,\big|\, \delta U \in \mathbb{R}^{k \times d},\;
	U \delta U^\top = 0,  \ X = U^\top Y \Bigr\},$$
	respectively.
	Using the decomposition of the tangent space of $ \mathcal{M}_k(\mathbb{R}^{d \times M})$ at the point $X=U^\top Y$, we can express the orthogonal projector onto this tangent space as the sum of the components in $V_X\mathcal{M}_k(\mathbb{R}^{d \times M})$, and $H_X\mathcal{M}_k(\mathbb{R}^{d \times M})$, namely $Z_V$ and $Z_H$, respectively.
	We write the orthogonal projector using the information
	onto $V_X \mathcal{M}_k(\mathbb{R}^{d \times M})$ and $H_X \mathcal{M}_k(\mathbb{R}^{d \times M})$
	with respect to the mixed scalar product
	\begin{equation}\label{eq: mixed scalar product}
		\langle A,B \rangle_{\mathrm{F},M} := \operatorname{Tr}\!\left( A^\top B \, \Sigma_M \right),
		\qquad \text{ for }A,B \in \mathbb{R}^{d \times M},
	\end{equation}
	where
	\begin{equation*}
		\Sigma_M = \operatorname{diag}\!\left( \frac{1}{M}, \ldots, \frac{1}{M} \right),
	\end{equation*}
	are the weights of the empirical scalar product. Relation \eqref{eq: mixed scalar product} translates into considering a weighted Frobenius scalar product
	defined through the weighted diagonal matrix whose elements are the weights of the Monte-Carlo discretization, acting on the corange of the application.
	
	We derive the projection onto $T_X \mathcal{M}_k(\mathbb{R}^{d \times M})$
	as minimization of the norm obtained by the scalar product \eqref{eq: mixed scalar product}. 
	To carry out the differentiation, we make repeated use of the identities \eqref{eq: matrix der}.
	For the projection onto the vertical component $Z_V$, we have
	\begin{equation*}
		\begin{aligned}
			Z_V &= \arg \min_{\delta Y} \| Z - U^\top \delta Y \|_{\mathrm{F},M} \\
			&= \arg \min_{\delta Y} \operatorname{Tr}\!\left( (Z - U^\top \delta Y)^\top (Z - U^\top \delta Y)\Sigma_M \right), \quad \text{for } Z \in \mathcal{M}_k(\mathbb{R}^{d \times M}).
		\end{aligned}
	\end{equation*}
	By taking the derivative with respect to the variation $\delta Y$, one obtains
	\begin{equation*}
		\frac{\partial}{\partial \delta Y} \| Z - U^\top \delta Y \|_{\mathrm{F},M}^2
		= - 2 U Z \Sigma_M + 2 \delta Y \Sigma_M = 0,
	\end{equation*}
	which implies that
	\begin{equation}\label{eq: deltaY}
		\delta Y = U Z.
	\end{equation}
	Therefore, the orthogonal projection onto the vertical space is given by
	\begin{equation}\label{eq: P_V}
		P_V[Z] = U^\top U Z, \quad \text{for } Z \in \mathcal{M}_k(\mathbb{R}^{d \times M}).
	\end{equation}
	
	Now we consider the projection onto the horizontal part $Z_H$, proceeding in the same fashion as in the vertical one one obtains:
	\begin{equation*}
		\begin{aligned}
			Z_H &= \arg \min_{\delta U^{\top}} \| Z - Z_V - \delta U^\top Y \|_{\mathrm{F},M} \\
			&= \arg \min_{\delta U^{\top}} \| Z - P_V Z  - \delta U^\top Y \|_{\mathrm{F},M} \\
			&= \arg \min_{\delta U^{\top}} \operatorname{Tr}\!\left(
			(P_U^\perp Z - \delta U^\top Y)^\top (P_U^\perp Z - \delta U^\top Y) \, \Sigma_M
			\right),  \quad \text{for } Z \in \mathcal{M}_k(\mathbb{R}^{d \times M}).
		\end{aligned}
	\end{equation*}
	Taking derivative with respect to $\delta U^\top$, one finds that
	\begin{equation*}
		\frac{\partial}{\partial \delta U^\top} \| Z - Z_V - \delta U^\top Y \|_{\mathrm{F},M}^2
		= 2 Y \Sigma_M Z^\top P_U^\perp - 2 Y \Sigma_M Y^\top \delta U = 0.
	\end{equation*}
	Thus,
	\begin{equation*}
		\delta U^\top = P_U^\perp Z \, \Sigma_M Y^\top \, \bigl( Y \Sigma_M Y^\top \bigr)^{-1}.
	\end{equation*}
    and, hence, the orthogonal projection onto the horizontal component reads as
	\begin{equation}\label{eq: P_H}
		P_H[Z] = P_U^\perp P_Y[Z],  \qquad \text{for } Z \in \mathcal{M}_k(\mathbb{R}^{d \times M}).
	\end{equation}
	where $P_Y$ is taken with respect to $\langle \ \cdot \ , \ \cdot \ \rangle_M$, i.e.\ the orthogonal projection onto the corange of $X$.
	
	Considering relations \eqref{eq: P_V} and \eqref{eq: P_H}, the full orthogonal projection onto $T_X \mathcal{M}_k(\mathbb{R}^{d \times M})$ is
	\begin{equation*}
		\begin{aligned}
		P_X[Z] &= P_V[Z] + P_H[Z] = P_U[Z] + P_U^\perp P_Y[Z] \\
	&= U^{\top} U Z + \left(I_{d \times d} - U^{\top}U \right) Z \Sigma_M Y^\top (Y \Sigma_M Y^\top)^{-1} Y.
		\end{aligned}
	\end{equation*}
	
	When one converts a SDE from the Stratonovich to the It\^o formulation, it is required to compute a correction term for the diffusion. As in standard projected SDEs on manifolds, this correction involves the differential of the projector $\mathrm{D}_X \big( P_X[Z] \big)$.
	Via chain rule and exploiting the fact that the projector $P_X$ is a linear operator, one has
	\begin{equation}\label{eq: deriv P_X}
	\mathrm{D}_X \big( P_X[Z] \big) = (\mathrm{D}_X P_X)[Z] + P_X[\mathrm{D}_X Z],
	\end{equation}
	where the last term on the right-hand side vanishes if $Z$ does not depend on the point $X$.
	
	For the sake of completeness, we compute the precise expression with respect to $\delta U$ and $\delta Y$ for this correction term. We have:
    \begin{equation}\label{eq: D_X P_X draft}
    	\begin{aligned}
	\mathrm{D}_X \left(P_U[Z] + P_U^\perp P_Y[Z]\right) &=  (\delta U^\top U + U^\top \delta U)[Z] + (0 - \delta U^\top U - U^\top \delta U) P_Y[Z] 
		+ P_U^\perp \, \delta P_Y[Z]  \\
		& = (\delta U^\top U + U^\top \delta U) P_Y^{\perp}[Z] 
		+ P_U^\perp \, \delta P_Y[Z]. \\
		\end{aligned}
	\end{equation}
		
	We want to make explicit the variations of the projection onto the corange $P_Y$. Via standard computations, we have
		\begin{equation*}\begin{aligned}
		\delta P_Y[Z] &= Z \Sigma_M (\delta Y)^\top (Y \Sigma_M Y^\top)^{-1} Y + Z \Sigma_M Y^\top \delta \left((Y \Sigma_M Y^\top)^{-1}\right) Y + Z \Sigma_M Y^\top (Y \Sigma_M Y^\top)^{-1} \delta Y.
		\end{aligned}\end{equation*}
	The variation with respect to the inverse of the Gramian is
		\begin{equation*}
	 \delta (Y \Sigma_M Y^\top)^{-1} = -(Y \Sigma_M Y^\top)^{-1} \, \delta (Y \Sigma_M Y^\top) \, (Y \Sigma_M Y^\top)^{-1},
	\end{equation*}
	and
	\begin{equation*}
	\delta (Y \Sigma_M Y^\top) = \delta Y \Sigma_M Y^\top + Y \Sigma_M \delta Y^\top.
	\end{equation*}
	Therefore, one has
	\begin{equation*}\begin{aligned}
	 \delta (Y \Sigma_M Y^\top)^{-1} 
		&= -(Y \Sigma_M Y^\top)^{-1} (\delta Y \Sigma_M Y^\top)(Y \Sigma_M Y^\top)^{-1} \\
		&\quad - (Y \Sigma_M Y^\top)^{-1} (Y \Sigma_M \delta Y^\top)(Y \Sigma_M Y^\top)^{-1}.
		\end{aligned}\end{equation*}
		Putting all the terms together, we obtain the following expression for the differential of the orthogonal projector onto the corange of $X$:
		\begin{equation*}
			\begin{aligned}
	\delta P_Y[Z] = & Z \Sigma_M (\delta Y)^\top (Y \Sigma_M Y^\top)^{-1} Y + Z \Sigma_M Y^\top (Y \Sigma_M Y^\top)^{-1} \delta Y\\
		& + Z \Sigma_M Y^\top \big( -(Y \Sigma_M Y^\top)^{-1} (\delta Y \Sigma_M Y^\top)(Y \Sigma_M Y^\top)^{-1} \big) Y \\
		&+ Z \Sigma_M Y^\top \big( -(Y \Sigma_M Y^\top)^{-1} (Y \Sigma_M \delta Y^\top)(Y \Sigma_M Y^\top)^{-1} \big) Y.
		\end{aligned}\end{equation*}
		After these considerations, \eqref{eq: D_X P_X draft} becomes
		\begin{equation}\label{eq: D_X P_X}
			\begin{aligned}
		D_X P_X[Z] &= (\delta U U^\top + U \delta U^\top) P_Y^\perp[Z] \\
		&\quad + P_U^\perp \Big[ Z \Sigma_M (\delta Y)^\top (Y \Sigma_M Y^\top)^{-1} Y \\
		&\qquad - Z \Sigma_M Y^\top (Y \Sigma_M Y^\top)^{-1} (\delta Y \Sigma_M Y^\top)(Y \Sigma_M Y^\top)^{-1} Y \\
		&\qquad - Z \Sigma_M Y^\top (Y \Sigma_M Y^\top)^{-1} (Y \Sigma_M \delta Y^\top)(Y \Sigma_M Y^\top)^{-1} Y \\
		&\qquad + Z \Sigma_M Y^\top (Y \Sigma_M Y^\top)^{-1} \delta Y \Big].
		\end{aligned}
	\end{equation}
		
		\subsection{The projected Stratonovich-SDE-based DLRA}
		As already discussed, we want to derive our DLRA equations via exploiting the chain rule of the Stratonovich integral. To compare the sought equations with the standard ones for Itô SDEs in order to see if they differ, we first start from a general Itô SDE, we convert it into the Stratonovich formulation, then we build our DLRA surrogate, and, finally, we convert it back to the Itô formalism.
		
		Our starting point is the SDE in Itô form\eqref{eq:SDE-int}, which we recall here for the sake of convenience,
		\begin{equation*}
			\begin{aligned}
			\mathrm{d}X_t &= a(t,X_t)\mathrm{d}t+b(t,X_t)\mathrm{d}W_t \\
			&= a(t,X_t)\mathrm{d}t+\sum_{q=1}^m b_q(t,X_t)\mathrm{d}W_t^q, 
			\end{aligned}
		\end{equation*}
		where $b_q(t,x)$ denotes the $q$-th column of $b(t,x)$. Assuming that $b_q$ is differentiable with respect to the spatial variable, the equivalent Stratonovich form is
		\begin{equation}\label{eq: stratonovich sde}
				\begin{aligned}
			\mathrm{d}X_t & = \widetilde{a}(t,X_t)\mathrm{d}t+b(t,X_t) \circ \mathrm{d}W_t \\
			&= \widetilde{a}(t,X_t)\mathrm{d}t+\sum_{q=1}^m b_q(t,X_t)\circ\mathrm{d}W_t^q 
		\end{aligned}
		\end{equation}
		where the Stratonovich drift $\widetilde{a}$ is
		\begin{equation*}
			\widetilde{a}(t,x)=a(t,x)-\frac{1}{2}\sum_{q=1}^m \partial_xb_q(t,x)\bigl[b_q(t,x)\bigr].
		\end{equation*}
		
		We consider $M$ particles $\left((X^{(M,i)}_t)_{t\geq 0}\right)_{i=1,\dots,M}$, where each $i$-th particle $(X^{(M,i)}_t)_{t\geq 0}$ is originated by the $i$-th realization of the same initial condition $X_0$ and is driven by a Brownian motion $W_t^{(i)}$, independent on all the others, i.e.\  $W_t^{(i)}\perp W_t^{(j)}$ for $i\neq j$, $i,j = \{1, \dots, M\}$. More specifically, for all $i=1,\dots, M$, each particle $(X^{(M,i)}_t)$ is characterized by the following properties:
		\begin{itemize}
			\item the initial condition of the SDE describing $X^{(M,i)}$ is $X_0^{(M,i)}= X_0(\omega_i)$ for the $i$-th realization $\omega_i \in \Omega$, where each $\omega_i $ is sampled independently from $\Omega$;
	        \item the component $\widetilde{a}\big(t, X^{(M,i)}_t\big)$ is the drift of the Stratonovich SDE \eqref{eq: stratonovich sde} computed in the $i$-th particle;
	       \item the component $b\big(t, X^{(M,i)}_t\big)$ is the diffusion of the Stratonovich SDE \eqref{eq: stratonovich sde} computed in the $i$-th particle;
	       \item $W_t^{(i)}$ is the $m$-dimensional Brownian Motion associated to the evolution of the $i$-th particle assumed in a Stratonovich form, where $(W_t^{(i)})_{i=1,\dots,M}$ is a set of $M$ independent Brownian motions;
	       \item $X_t^{(1,\ldots,M)}$ is the matrix whose 
	       $i$-th column is the $i$-th particle $X_t^{(M,i)}$, i.e.\
	       $$X_t^{(1,\ldots,M)}:= \left[X_t^{(M,1)} X^{(M,2)}_t \dots X^{(M,M)}_t \right] \in \mathbb{R}^{d \times M}.$$ Then, each particles is constrained to remain into the manifold $\mathcal{M}_k(\mathbb{R}^{d \times M})$ via projecting each drift $\widetilde{a}\big(t, X^{(M,i)}_t\big)$ and diffusion $b\big(t, X^{(M,i)}_t\big)$ into the tangent space of $ \mathcal{M}_k(\mathbb{R}^{d \times M})$ computed in the point $X_t^{(1,\ldots,M)}$ that resembles the set of all the particles $(X^{(M,i)}_t)$.
	\end{itemize}
		
	In the view of these considerations, the evolution of each particle can be described by the following matrix Stratonovich SDE on the manifold $\mathcal{M}_k(\mathbb{R}^{d \times M})$ \cite{hsu2002stochastic}:
	\begin{equation}\label{eq: stratonovic DLRA}
		\begin{aligned}
		 \mathrm{d}X_t^{(M,i)} &= P_{X_t^{(1,\ldots,M)}} \bigl[\widetilde{a}\bigl(t, X_t^{(M,i)}\bigr) e_i^{\top}\bigr] \mathrm{d}t 
		+ P_{X_t^{(1,\ldots,M)}} b\bigl(t, X_t^{(M,i)}\bigr) \circ \mathrm{d}W_t^{(i)} \\
		& =  P_{X_t^{(1,\ldots,M)}}\bigl[\widetilde a(t,X_t^{(M,i)})e_i^{\top}\bigr]\mathrm{d}t+\sum_{q=1}^mP_{X_t^{(1,\ldots,M)}}\bigl[b_q(t,X_t^{(M,i)})e_i^\top\bigr]\circ\mathrm{d}W_t^{(i,q)}, \qquad i=1,\dots,M,
		\end{aligned}
	\end{equation}
	where $\circ$ denotes the stochastic integral in a Stratonovich form, $P_{X_t^{(1,\ldots,M)}}$ is the orthogonal projection onto the tangent space 
	$T_{X_t^{(1,\ldots,M)}} \mathcal{M}_k(\mathbb{R}^{d \times M})$, $e_i^{\top}$ is the canonical basis in $\mathbb{R}^M$, and $W_t^{(i,q)}$ is the $q$-th coordinate of the $i$-th Brownian motion. Relation \eqref{eq: stratonovic DLRA} is saying that the evolution of all the particles is constrained to belong to $\mathcal{M}_k(\mathbb{R}^{d \times M})$ through the projection of the dynamics onto $T_{X_t^{(1,\ldots,M)}} \mathcal{M}_k(\mathbb{R}^{d \times M})$.  For the sake of notation, we define	
	\begin{equation*}
		\tilde{A}_i(t,X_t)= \tilde{a}(t,X_t^{(M,i)})e_i^{\top}, \quad B_{i,q}(t,X_t):=b_q(t,X_t^{(M,i)})e_i^\top\in\mathbb{R}^{d\times M},
	\end{equation*}
	
	One can define  $P_{X_t^{(1,\ldots,M)}}$ following priors computations provided in Section \ref{sec: orth proj}. Indeed, via using the previous notation, we can write
	\begin{equation*}
		P_{X_t^{(1,\ldots,M)}}[Z] = P_{U_t(X_t^{(1,\ldots,M)})}^\perp P_{Y_t(X_t^{(1,\ldots,M)})}[Z] 
		+ P_{U_t(X_t^{(1,\ldots,M)})}[Z], \quad \text{ for } Z \in \mathcal{M}_k(\mathbb{R}^{d \times M}),
	\end{equation*}
	where $P_{U_t(X_t^{(1,\ldots,M)})}$ and $P_{Y_t(X_t^{(1,\ldots,M)})}$ are the orthogonal projectors onto the range and the corange of $X_t^{(1,\ldots,M)}$, respectively.	
	For the sake of notation, we denote
	\begin{equation*}
		P_{X_t^M} = P_{X_t^{(1,\ldots,M)}}, \quad 
		P_{U_t^M} = P_{U_t(X_t^{(1,\ldots,M)})}, \quad 
		P_{Y_t^M} = P_{Y_t(X_t^{(1,\ldots,M)})} .
	\end{equation*}
	and, hence, then equation \eqref{eq: stratonovic DLRA} can be rewritten as
	\begin{equation*}
		\begin{aligned}
		 \mathrm{d}X_t^{(M,i)} & = P_{X_t^M} \bigl[\tilde{A}_i(t,X_t)\bigr] \mathrm{d}t 
		+ P_{X_t^M} \bigl[b\bigl(t, X_t^{(M,i)}\bigr) e_i^{\top} \bigr]\circ \mathrm{d}W_t^{(i)} \\
		& = P_{X_t^M}\bigl[\tilde{A}_i(t,X_t)\bigr]\mathrm{d}t+\sum_{q=1}^mP_{X_t^M}\bigl[B_{i,q}(t,X_t)\bigr]\circ\mathrm{d}W_t^{(i,q)}, \qquad i=1,\dots,M.
		\end{aligned}
	\end{equation*}
	
	Assuming the derivability with respect to the spatial variable of the diffusion $b$, then the SDE \eqref{eq: stratonovic DLRA} written in a Stratonovich form can be translated into an equivalent Itô formulation as follows \cite{oksendal2003stochastic}:
	\begin{equation}\label{eq: ito DLRA}
		\begin{aligned}
			 \mathrm{d}X_t^{(M,i)} &= \left(P_{X_t} \bigl[\tilde{A}_i(t,X_t)\bigr]+ \frac{1}{2} \sum_{q=1}^m D_{X_t^{(M,i)}} \bigl[ P_{X_t^M} B_{i,q}(t,X_t) \bigr] : \bigl[P_{X_t^M} B_{i,q}(t,X_t)\bigr] \, \right) \mathrm{d}t \\
			&\quad + P_{X_t^M} \bigl[b\bigl(t, X_t^{(M,i)}\bigr) e_i^{\top}\bigr] \mathrm{d}W_t^{(i)}, \quad \text{ for all } i=1, \dots, M.
		\end{aligned}
	\end{equation}
	
	Therefore, the correction term from the Stratonovich to the Itô form reads as
	\begin{equation*}
	\frac{1}{2} \sum_{q=1}^m \mathrm{D}_{X_t^{(M,i)}} \bigl[ P_{X_t^M} B_{i,q}(t,X_t) \bigr] : \bigl[P_{X_t^M} B_{i,q}(t,X_t)\bigr].
	\end{equation*}
	Making explicit $\widetilde{a}$ and the derivatives in the correction term \eqref{eq: deriv P_X}, one obtains 
	\begin{equation}\label{eq: i-th part DLRA}
		\begin{aligned}
			\mathrm{d}X_t^{(M,i)} 
			= & \Bigg(P_{X_t^M} \, \bigl[A\bigl(t, X_t^{(M,i)}\bigr) \bigr]
			- \frac{1}{2} \sum_{q=1}^m \left[P_{X_t^M} \bigl[ \mathrm{D}_{\mathrm{D}_{X_t^{(M,i)}}} B_{i,q}(t,X_t) \bigr] \right] : \bigl[B_{i,q}(t,X_t)\bigr] \\ 
			&+	\frac{1}{2} \sum_{q=1}^m D_{X_t^{(M,i)}} \bigl[ P_{X_t^M} B_{i,q}(t,X_t) \bigr] : \bigl[P_{X_t^M} B_{i,q}(t,X_t)\bigr] \, \Bigg) \mathrm{d}t+ P_{X_t^M} \, b\bigl(t, X_t^{(M,i)}\bigr) \mathrm{d}W_t^{(i)}\\
				= & \Bigg(P_{X_t^M} \, a\bigl(t, X_t^{(M,i)}\bigr) 
			- \frac{1}{2} \sum_{q=1}^m \left[P_{X_t^M} \bigl[ \mathrm{D}_X B_{i,q}(t,X_t) \bigr] \right] : \bigl[P_{X_t^M}^{\perp}B_{i,q}(t,X_t)\bigr] \\ 
			&+\frac{1}{2} \sum_{q=1}^m \left[ \mathrm{D}_X P_{X_t^M} \bigl[ P_{X_t^M} B_{i,q}(t,X_t)\bigr] \bigl[ B_{i,q}(t,X_t) \bigr] \right]  \, \Bigg) \mathrm{d}t+ P_{X_t^M} \, b\bigl(t, X_t^{(M,i)}\bigr) \mathrm{d}W_t^{(i)}, \quad i=1,\dots,M,
		\end{aligned}
	\end{equation}
	where we define $A_i(t,X_t)= a(t,X_t^{(M,i)})e_i^{\top}$.
	
	\begin{Remark}
	In the view of \eqref{eq: i-th part DLRA}, we can see the whole particle system as a matrix Itô SDE:
	\begin{equation*}
		\begin{aligned}
			\mathrm{d}X_t^{(1,\dots,M)} 
			=& P_{X_t^{(1,\dots,M)}} 
			\begin{bmatrix}
				a(t, X_t^{(M,1)}) & \dots & a(t, X_t^{(M,M)})
			\end{bmatrix} \mathrm{d}t\\
			& - \frac{1}{2} \sum_{q=1}^m \sum_{\ell=1}^d 
			P_{X_t^{(1,\dots,M)}} 
			\begin{bmatrix}
				\mathrm{d}_{x_\ell} B_{1,q}(t, X_t) & \dots & \mathrm{d}_{x_\ell} B_{M,q}(t, X_t))
			\end{bmatrix}\\
			& \qquad \qquad \odot
			\begin{bmatrix} 
				P_{X_t^{(1,\dots,M)}}^{\perp} 
			\begin{bmatrix}
				 B_{1,q}(t, X_t) & \dots & B_{M,q}(t, X_t))
				\end{bmatrix} \end{bmatrix} \mathrm{d}t \\
			& + \frac{1}{2} \sum_{q=1}^m \sum_{\ell=1}^d 
			[D_X P_{X_t^{(1,\dots,M)}} ]
			\begin{bmatrix}
				P_{X_t^{(1,\dots,M)}} B_{1,q}(t, X_t) & \dots & P_{X_t^{(1,\dots,M)}} B_{M,q}(t, X_t))
			\end{bmatrix}\\
			& \qquad \qquad \odot
				\begin{bmatrix}
					B_{1,q}(t, X_t) & \dots & B_{M,q}(t, X_t))
			\end{bmatrix} \mathrm{d}t \\
			& + \sum_{q=1}^m 
			P_{X_t^{(1,\dots,M)}} 
			\begin{bmatrix}				
				B_{1,q}(t, X_t) & \dots & B_{M,q}(t, X_t))
			\end{bmatrix}
			\odot
			\begin{bmatrix}
				\mathrm{d}W_t^{(1,q)}  \dots \mathrm{d}W_t^{(M,q)}
			\end{bmatrix},
		\end{aligned}
	\end{equation*}
	where we recall that $B_{i,q}(t, X_t)$ is the $q$-th column of the diffusion computed in the $i$-th particle $X_t^{(M,i)}$, $\odot$ denotes the Hadamard product, $W_t^{(i),q}$ is the $q$-th coordinate of the $i$-th Brownian motion, and the following relations hold
	\begin{equation*}
		P_{X_t^{(1,\dots,M)}} : \mathbb{R}^{d \times M} \to \mathbb{R}^{d \times M},
	\end{equation*}
	\begin{equation*}
			\begin{bmatrix}
			B_{1,q}(t, X_t) & \dots & B_{M,q}(t, X_t))
		\end{bmatrix} = \begin{bmatrix}
			b_q(t, X_t^{(M,1)}) & \dots & b_q(t, X_t^{(M,M)})
		\end{bmatrix} \in \mathbb{R}^{d \times M}, 	\quad \forall q=1,\dots,m,
	\end{equation*}
	\begin{equation*}
		\bigl[(b(t, X_t^{(M,i)}))_{i=1,\dots,M}\bigr] =
		\begin{bmatrix}
			b(t, X_t^{(M,1)}) & \dots & b(t, X_t^{(M,M)})
		\end{bmatrix} \in \mathbb{R}^{d \times m \times M}, 
	\end{equation*}
	\begin{equation*}
		\bigl[ \mathrm{D}_X b(t, X_t^{(M,i)}) \bigr] \in \mathbb{R}^{d \times d \times m},
		\quad 	\bigl[ \mathrm{d}_{x_\ell} b(t, X_t^{(M,i)}) \bigr] \in \mathbb{R}^{d \times d}
	\end{equation*}
	\begin{equation*}
		\bigl[ \mathrm{d}_{x_\ell} B_{i,q}(t, X_t)  \bigr] \in \mathbb{R}^d, \quad \forall q=1,\dots,m, \quad \quad \forall i=1,\dots,M,
	\end{equation*}
	where $\mathrm{d}_{x_\ell}$ denotes the differential with respect to the $\ell$ physical component, i.e.\ 
	$$\mathrm{d}_{x_\ell} b(t, X_t^{(M,i)}) = [\mathrm{d}_{x_\ell}  b_1(t, X_t^{(M,i)}) \ \mathrm{d}_{x_\ell} b_2(t, X_t^{(M,i)}) \  \dots \ \mathrm{d}_{x_\ell} b_m(t, X_t^{(M,i)}) ]^{\top} \in \mathbb{R}^{m}.$$
    \end{Remark}
	
		One can give an explicit expression to the additional term coming from the Itô to Stratonovich correction 
	\begin{equation*}
		\begin{aligned}
			D_{X_t^{(M,i)}}\left(P_{X_t^{M}}\bigl[B_{i,q}(t,X)\bigr]\right)\left[P_{X_t^{M}}\bigl[B_{i,q}(t,X)\bigr]\right] 
			&= D_{X_t^{(M,i)}} P_{X_t^{M}} \left[ P_{X_t^{M}}\bigl[B_{i,q}(t,X)\bigr]\right]\bigl[B_{i,q}(t,X)\bigr] \\
			&\quad + P_{X_t^{M}} \left[D_{X_t^{(M,i)}}  B_{i,q}(t,X)\left[P_{X_t^{M}} \bigl[B_{i,q}(t,X)\bigr]\right]\right].
		\end{aligned}
	\end{equation*}
	via denoting the vertical and horizontal components coordinate of  $B_{i,q}(t,X)$.
	Since $B_{i,q}(t,X)$ depends only on the $i$-th column of $X^{(1,\ldots,M)}$, the second term is localized on the $i$-th particle. 
	Thus, for each projected noise direction
	\begin{equation*}
		Z_{i,q}:=P_{X_t^{M}}\bigl[B_{i,q}(t,X)\bigr] \in \mathbb{R}^{d \times m},
	\end{equation*}
	the corresponding vertical and horizontal component factors $\delta_{i,q}Y_t$ and $\delta_{i,q}U_t$ are, respectively,
	\begin{equation}\label{eq: horizontal component factor}
		\begin{aligned}
				\delta_{i,q}Y_t & = G^Y_{i,q} = U^M_t Z_{i,q}=U^M_t B_{i,q}(t,X) \in  \mathbb{R}^{k \times m} \\
			\delta_{i,q}(U_t^M)^\top  &=(G^U_{i,q})^\top = P_{U_t^M}^\perp Z_{i,q}\Sigma_M(Y^M_t)^\top\bigl(Y^M_t\Sigma_M(Y^M_t)^{\top}\bigr)^{-1}\\
			&=P_{U_t^M}^\perp B_{i,q}(t,X)\Sigma_M(Y^M_t)^{\top}\bigl(Y^M_t\Sigma_M(Y^M_t)^{\top}\bigr)^{-1} \in \mathbb{R}^{d \times m \times k},
		\end{aligned}
	\end{equation}
	for all $i=1,\dots, M$.
	Therefore, for each additional component on the drift we can write
	\begin{equation}\label{eq: D_X P first term}
		\begin{aligned}
			&D_XP_{X_t^M}\bigl[P_{X_t^M}\bigl[B_{i,q}(t,X)\bigr]\bigr]\bigl[B_{i,q}(t,X_t)\bigr] \\
			=& \bigl(\delta_{i,q}(U_t^M)^\top U_t^M+(U_t^M)^\top\delta_{i,q}U_t^M\bigr)P_{Y_t^M}^\perp\bigl[B_{i,q}(t,X_t)\bigr] \\
			& +P_{U_t^M}^\perp\Bigl[B_{i,q}(t,X_t)\Sigma_M(\delta_{i,q}Y_t^M)^\top\bigl(Y_t^M\Sigma_M(Y_t^M)^\top\bigr)^{-1}Y_t^M \\
			&\quad -B_{i,q}(t,X_t)\Sigma_M(Y_t^M)^\top\bigl(Y_t^M\Sigma_M(Y_t^M)^\top\bigr)^{-1}
			\bigl(\delta_{i,q}Y_t^M\Sigma_M(Y_t^M)^\top\bigr)\bigl(Y_t^M\Sigma_M(Y_t^M)^\top\bigr)^{-1}Y_t^M \\
			&\quad -B_{i,q}(t,X_t)\Sigma_M(Y_t^M)^\top\bigl(Y_t^M\Sigma_M(Y_t^M)^\top\bigr)^{-1}
			\bigl(Y_t^M\Sigma_M(\delta_{i,q}Y_t^M)^\top\bigr)\bigl(Y_t^M\Sigma_M(Y_t^M)^\top\bigr)^{-1}Y_t^M \\
			&\quad +B_{i,q}(t,X_t)\Sigma_M(Y_t^M)^\top\bigl(Y_t^M\Sigma_M(Y_t^M)^\top\bigr)^{-1}\delta_{i,q}Y_t^M\Bigr]. \\
			=& P_{U_t^M}^\perp B_{i,q}(t,X)\Sigma_M(Y^M_t)^\top\bigl(Y^M_t\Sigma_M(Y^M_t)^\top\bigr)^{-1} U_t^M P_{Y_t^M}^\perp\bigl[B_{i,q}(t,X_t)\bigr] \\
			&+(U_t^M)^{\top}\bigl(Y^M_t\Sigma_M(Y^M_t)^{\top}\bigr)^{-1} Y^M_t \Sigma_M B_{i,q}(t,X)^{\top} P_{U^M_t}^\perp \bigr)P_{Y_t^M}^\perp\bigl[B_{i,q}(t,X_t)\bigr] \\
			& +P_{U_t^M}^\perp\Bigl[B_{i,q}(t,X_t)\Sigma_M B_{i,q}(t,X)^\top (U^M_t)^{\top}\bigl(Y_t^M\Sigma_M(Y_t^M)^\top\bigr)^{-1}Y_t^M \\
			&\quad -B_{i,q}(t,X_t)\Sigma_M(Y_t^M)^\top\bigl(Y_t^M\Sigma_M(Y_t^M)^\top\bigr)^{-1}
			\bigl(U^M_t B_{i,q}(t,X)\Sigma_M(Y_t^M)^\top\bigr)\bigl(Y_t^M\Sigma_M(Y_t^M)^\top\bigr)^{-1}Y_t^M \\
			&\quad -B_{i,q}(t,X_t)\Sigma_M(Y_t^M)^\top\bigl(Y_t^M\Sigma_M(Y_t^M)^\top\bigr)^{-1}
			\bigl(Y_t^M \Sigma_M  B_{i,q}(t,X)^\top (U^M_t)^{\top}\bigr)\bigl(Y_t^M\Sigma_M(Y_t^M)^\top\bigr)^{-1}Y_t^M \\
			&\quad +B_{i,q}(t,X_t)\Sigma_M(Y_t^M)^\top\bigl(Y_t^M\Sigma_M(Y_t^M)^\top\bigr)^{-1}U^M_t B_{i,q}(t,X)\Bigr]. \\
		\end{aligned}
	\end{equation}
	
	\subsection{The equations for particle-based $U$ and $Y$}
		
	Now, one wants to find reasonable equations for the components of the deterministic basis $U^M$ and the stochastic one $Y^M$ of a possible DO solution such that $X^M_t = (U^M_t)^{\top} Y^M_t$. Then, we would formally take the limit for $M \to \infty$ in order to obtain differential equations for the deterministic basis $U$ and the stochastic one $Y$ for standard SDEs. Would this search be successful, one would obtain equations that can be computable and, hence, their solution be assembled to obtain a DLRA $X$. For the sake of notation, we omit the superscript $M$ unless further precision is needed. 
	
	 In contrast with \cite{kazashi2025dynamical}, we do not restrict $U^M$ to follows a deterministic ODE, but we allow to have a diffusion term in addition to the deterministic drift. By ansatz, we ask for $U$ and $Y$ to be described by the following SDEs:
	\begin{equation*}
	 \mathrm{d}Y^M_t = \alpha_Y \mathrm{d}t + \beta_Y \mathrm{d}W_t, \quad 	\mathrm{d}U^M_t = \alpha_U \mathrm{d}t + \beta_U \mathrm{d}W_t, 
	\end{equation*}
	where $U_t^M (U_t^M)^{\top} = I_{k\times k}$, i.e.\ $U_t^M$ has orthonormal rows, $U_t^M (\mathrm{d}U_t^M)^{\top} = 0$, with coefficients $\alpha_{U}\in\mathbb{R}^{k \times d}$ and $\beta_{U}\in\mathbb{R}^{k\times d \times m}$,  $\alpha_{Y}\in\mathbb{R}^{k}$ and $\beta_{Y}\in\mathbb{R}^{k\times m}$,
	where $\alpha_{U}$ is deterministic, $\beta_{U}$, $\alpha_{Y}$ and $\beta_{Y}$ are 
	progressively measurable
	and have continuous paths almost surely. As usual, we ask for the linear independence of the components of $Y^{M}_t$ for any $t \in [0,T]$, i.e.\ $Y^M_t\Sigma_M(Y^M_t)^\top$ is invertible and, hence of rank equal to $k$.
	
By the Itô formula associated to  $X^M_t = (U^M_t)^{\top} Y^M_t$ and via relation \eqref{eq: horizontal component factor} we already obtain relations for the diffusion $\beta_{Y}$ and $\beta_{U}$
	\begin{equation*}
		\begin{aligned}
			\beta_{Y} &=\sum_{i=1}^M\sum_{q=1}^m G^Y_{i,q}(t,X_t) = \sum_{i=1}^M\sum_{q=1}^m U B_{i,q}(t,X),\\
			\beta_{U} &= \sum_{i=1}^M\sum_{q=1}^m G^U_{i,q}(t,X_t) = \sum_{i=1}^M\sum_{q=1}^m P_{U_t^M}^\perp B_{i,q}(t,X)\Sigma_M(Y^M_t)^{\top}\bigl(Y^M_t\Sigma_M(Y^M_t)^{\top}\bigr)^{-1}.
		\end{aligned}
	\end{equation*}
	On the other hand, to obtain relations for $\alpha_Y$ and $ \alpha_U$, we again exploit the Itô formula for the $i$-th particle
	\begin{equation}\label{eq: ito prod}
		\mathrm{d}X^{(M,i)}_t=(\mathrm{d}U^M_t)^\top Y^{(M,i)}_t+(U_t^M)^\top \mathrm{d}Y^{(M,i)}_t+\mathrm{d}\langle (U^M)^\top,Y^{(M,i)}\rangle_t.
	\end{equation}
	Via multiplying by $U^M_t$ both members of the last relation and exploiting the gauge condition $U^M,\mathrm{d}(U^M)^\top = 0$ and the orthogonality of the row of $U^M_t$, we obtain the following relation for $i$-th particle of $\mathrm{d}Y^M_t$
	\begin{equation}\label{eq: Y_M}
		\begin{aligned}
		\mathrm{d}Y^{(M,i)} _t = & \Bigg(U^M_t a(t, X_t^{(M,i)} ) + \sum_{q=1}^m \bigl(Y_t^M\Sigma_M(Y_t^M)^\top\bigr)^{-1} Y_t^M \Sigma_M b_{q}(t,X)^{\top} P_{U^M_t}^\perp \bigr)P_{Y_t^M}^\perp\bigl[B_{i,q}(t,X_t)\bigr] \\
		&-  \frac{1}{2} \sum_{q=1}^m U^M_t [ D_XB_{i,q}(t, X_t)][ P_{X_t^{(1,\dots,M)}}^{\perp} B_{i,q}(t, X_t)]   \Bigg) \mathrm{d}t + P_{U^M_t} b(t, X_t^{(M,i)} ) \mathrm{d} W_t^{i}
		\end{aligned}
	\end{equation}
	
	Consequently, the equations for $U^M$ is derived inserting \eqref{eq: Y_M} in \eqref{eq: ito prod}, multiplying for $\Sigma_M Y_t^M$, and assuming the Gramian invertible. Then, one has
	\begin{equation*}
		\begin{aligned}
			\mathrm{d}U^M_t =& \bigl(Y_t^M\Sigma_M(Y_t^M)^\top\bigr)^{-1} \Bigg( Y_t^M\Sigma_M a(t, X_t)^\top  -  \frac{1}{2} \sum_{q=1}^m [( Y_t^M \Sigma_M b_q(t, X_t)^\top)]  [P_{Y^M_t}[ D_Xb_{q}(t, X_t)]]^{\top} \\
			&-  \frac{1}{2} \sum_{q=1}^m U^M_t b_{q}(t,X) \Sigma_M (Y^M_t)^{\top} \bigl(Y_t^M\Sigma_M(Y_t^M)^\top\bigr)^{-1}  Y_t^M \Sigma_M b_{q}(t,X)^{\top} \Bigg) P_{U^M_t}^{\perp} \mathrm{d}t \\
			& + \frac{1}{2} \sum_{q=1}^m  \Bigg( U_t^M b_q(t,X_t) \Sigma_M b_q(t,X_t)^{\top}  \\
			& \qquad - \Big( Y_t^M \Sigma_M b_q(t, X_t)^{\top} (U^M_t)^{\top} + U_t^M b_q(t,X_t) \Sigma_M Y_t^M\Big) (Y_t^M \Sigma_M (Y_t^M)^{\top})^{-1} Y_t^M \Sigma_M b_q(t, X_t)^{\top} \\
			&\qquad +  Y_t^M \Sigma_M b_q(t,X_t)^{\top} (U^M_t)^{\top} (Y_t^M \Sigma_M (Y_t^M)^{\top})^{-1}  Y_t^M \Sigma_M b_q(t,X_t)^{\top} \Bigg) P_{U^M_t}^\perp \mathrm{d}t \\
			&+ \bigl(Y_t^M\Sigma_M(Y_t^M)^\top\bigr)^{-1} \sum_{i=1}^M\sum_{q=1}^m Y_t^M\Sigma_MB_{i,q}(t,X_t)^\top P_{U_t^M}^\perp\mathrm{d}W_t^{(i,q)}.
		\end{aligned}
	\end{equation*}

	\subsection{Monte-Carlo convergence of the particle system}
	
	We now pass formally to the limit $M\to+\infty$. In order to compute this limit, we assume by ansatz that this limit exists and convergence of this limit holds. A possible setting of assumptions that guarantees these hypothesis is the one of \cite{kazashi2026dynamicalpartII}.
	
	Formally, we assume that $Y_t^{(M,i)}$ converges in $L^2(\Omega)$ to the mean-field limit process $Y_t$ for each particle $i$ and the same holds for $U^M_t$, which converges to the basis $U_t$, which will be proved to be deterministic. We want to derive equations for this limits and retrieve the limit DLRA $X$ via Itô formula.
	
	Recall that $Y^M_t\Sigma_M(Y^M_t)^{\top}=\frac{1}{M}Y^M_t(Y^M_t)^\top..$ Then, via Monte Carlo approximation we suppose that this empirical Gramian of the stochastic basis converges to one of the continuous-in-time process, i.e.\
	\begin{equation*}
		Y_t^M\Sigma_M(Y_t^M)^\top \underset{M\to+\infty}{\longrightarrow} C_{Y_t}:=\mathbb{E}\bigl[Y_tY^\top_t\bigr],
	\end{equation*}
	Similarly, the projector onto the corange becomes the $L^2$-projection
	\begin{equation*}
		P_{Y_t^M}[Z]  \underset{M\to+\infty}{\longrightarrow} P_{Y_t}[Z] =\mathbb{E}\bigl[ZY_t^\top\bigr]C^{-1}_{Y_t} Y_t,
	\end{equation*}
	and, hence, the limiting tangent projector on the rank-$k$ manifold in $L^2(\Omega,\mathbb{R}^d)$ is
	\begin{equation*}
	P_{X_t}[Z]=P_{U_t}[Z]+P_{U_t}^\perp\mathbb{E}\bigl[ZY^\top_t\bigr]C^{-1}_{Y_t} Y_t.
	\end{equation*}
	
	For the noise terms, observe that
	\begin{equation}\label{eq: noise U^M}
		\begin{aligned}
		\sum_{q=1}^m \sum_{i=1}^{M} G^U_{i,q} = & \sum_{q=1}^m \sum_{i=1}^{M}  C_{Y^M_t}^{-1}Y^{(M,i)}_t\frac1M B_{i,q}^\top P_{U^M_t}^\perp\\
		=& \sum_{q=1}^mC_{Y^M_t}^{-1} \left( \frac{1}{M} \sum_{i=1}^{M} Y^{(M,i)}_t b_q\bigl(t,X^{(M,i)}_t\bigr)^\top \right) P_{U^M_t}^\perp \\
		\end{aligned}
	\end{equation}
	We consider a uniform partition of the time interval $[0,T]$, namely $\Delta := \{n \Delta t, \text{ with } \Delta t = \frac{T}{N}\}$. Suppose further that $t= n_{t} \Delta t$. We want to consider the limit for $M$ which goes to $+\infty$ for the Itô integral with integrand \eqref{eq: noise U^M}, i.e.\ 
	\begin{equation}\label{eq: noise U^M lim M}
		\begin{aligned}
			\lim\limits_{M \to \infty}\sum_{q=1}^m \sum_{i=1}^{M} \int_{0}^{t}G^U_{i,q} \mathrm{d}W_t^{(i,q)} = & \lim\limits_{M \to \infty} \sum_{q=1}^m \sum_{i=1}^{M} \int_{0}^{t} C_{Y^M_s}^{-1}Y^{(M,i)}_s\frac1M B_{i,q}^\top P_{U^M_s}^\perp \mathrm{d}W_s^{i,q}\\
			=&	\lim\limits_{M \to \infty} \sum_{q=1}^m  \int_{0}^{t} C_{Y^M_s}^{-1} \left( \frac{1}{M} \sum_{i=1}^{M} Y^{(M,i)}_s b_q\bigl(s,X^{(M,i)}_s\bigr)^\top \right) P_{U^M_s}^\perp \mathrm{d}W_s^{(i,q)} \\
			=&	\lim\limits_{M \to \infty} \sum_{q=1}^m  \lim\limits_{\substack{ \Delta t \to 0 \\ n_t \Delta t = t}}\sum_{n=0}^{n_t-1} C_{Y^M_{n \Delta t}}^{-1} \left( \frac{1}{M} \sum_{i=1}^{M} Y^{(M,i)}_{n \Delta t} b_q\bigl(n \Delta t,X^{(M,i)}_{n \Delta t}\bigr)^\top \right) P_{U^M_{n \Delta t}}^\perp \Delta W_n^{(i,q)}, \\
		\end{aligned}
	\end{equation}
	where $\{ \Delta W_n^{(i,q)}\}_n$ are all independent Brownian increments, i.e.\ $\Delta W_n^{(i,q)} \sim \mathcal{N}(0, \Delta t)$ for all $i=\{1,\dots, M\}$, $\forall q \in \{1,\dots,m\}$, $\mathbb{E}[\Delta W_n^{(i,q)}\Delta W_h^{(j,r)}]$ for all $n \neq h$, or $\forall i \neq j$, or $\forall q \neq r$.
	
	Then, we can write informally that
	\begin{equation}\label{eq: lim Y}
		\begin{aligned}
			Y_t(\omega_i) = &\lim\limits_{\Delta t \to 0} Y^{(M,i)}_{n \Delta t} \\
			= &\lim\limits_{\Delta t \to 0}  \left( Y^{(M,i)}_{0} + \sum_{h=0}^{n_t-1} \left( h_{Y^M_{\mathrm{drift}}}(h \Delta t, X^{(M,i)}_{h \Delta t}) \Delta t + \sum_{r=1}^m P_{U^M_{h \Delta t}} b_r(h \Delta t, X^{(M,i)}_{h \Delta t})  \Delta W_{h}^{(i,r)} \right) \right),
		\end{aligned}
	\end{equation}
	where $h_{Y^M_{\mathrm{drift}}}$ denotes the drift of $Y^{(M,i)}$.
	Plugging \eqref{eq: lim Y} in \eqref{eq: noise U^M lim M}, we find that
	\begin{equation*}
		\begin{aligned}
				&\lim\limits_{M \to \infty}\sum_{q=1}^m \sum_{i=1}^{M} \int_{0}^{t}G^U_{i,q} \mathrm{d}W_t^{(i,q)} \\
				=&	\lim\limits_{M \to \infty} \sum_{q=1}^m  \lim\limits_{\substack{ \Delta t \to 0 \\ n_t \Delta t = t}}\sum_{n=0}^{n_t-1} \Bigg( \frac{1}{M} \sum_{i=1}^{M} C_{Y^M_{n \Delta t}}^{-1}   \left( Y^{(M,i)}_{0} + \sum_{h=0}^{n_t-1} \left( h_{Y^M_{\mathrm{drift}}}(h \Delta t, X^{(M,i)}_{h \Delta t}) \Delta t + \sum_{r=1}^m U^M_{h \Delta t} b_r(h \Delta t, X^{(M,i)}_{h \Delta t})  \Delta W_{h}^{(i,r)} \right) \right) \\
				& \qquad \qquad \qquad \qquad  \qquad \qquad \qquad  \qquad \cdot b_q\bigl(n \Delta t,X^{(M,i)}_{n \Delta t}\bigr)^\top \Bigg) P_{U^M_{n \Delta t}}^\perp \Delta W_n^{(i,q)}, \\
		\end{aligned}
\end{equation*}
   
   Via independence of the Brownian increment, for the limit for $M \to \infty$, the only terms in the sum in $i$ that are no null are the one of the type
   \begin{equation}\label{eq: brownian incerement}
   	\mathbb{E}[ \Delta W_n^{(i,q)} \Delta W_n^{(i,q)}], \qquad \forall n \in \{1,\dots,N\}, \ \forall i \in \{1,\dots,M\}, \ \forall q \in \{1,\dots,m\}.
   \end{equation}
	Therefore, one has that
		\begin{equation*}
		\begin{aligned}
			&\lim\limits_{M \to \infty}\sum_{q=1}^m \sum_{i=1}^{M} \int_{0}^{t}G^U_{i,q} \mathrm{d}W_t^{(i,q)} \\
			=&	\lim\limits_{M \to \infty}  \lim\limits_{\substack{ \Delta t \to 0 \\ n_t \Delta t = t}}\sum_{n=0}^{n_t-1} \Bigg( \frac{1}{M} \sum_{i=1}^{M} C_{Y^M_{n \Delta t}}^{-1}   \left(\sum_{q=1}^m U^M_{n \Delta t} b_q(n \Delta t, X^{(M,i)}_{n \Delta t}) \Delta W_n^{(i,q)} \right) b_q\bigl(n \Delta t,X^{(M,i)}_{n \Delta t}\bigr)^\top \Bigg) P_{U^M_{n \Delta t}}^\perp \Delta W_n^{(i,q)}, \\
			=&	 \lim\limits_{\substack{ \Delta t \to 0 \\ n_t \Delta t = t}}\sum_{n=0}^{n_t-1} \lim\limits_{M \to \infty} C_{Y^M_t}^{-1} \Bigg( \frac{1}{M} \sum_{i=1}^{M}  \left(\sum_{q=1}^m U^M_{n \Delta t} b_q(n \Delta t, X^{(M,i)}_{n \Delta t}) \Delta W_n^{(i,q)} \right) b_q\bigl(n \Delta t,X^{(M,i)}_{n \Delta t}\bigr)^\top \Bigg) P_{U^M_{n \Delta t}}^\perp \Delta W_n^{(i,q)}, \\
			=&	 \lim\limits_{\substack{ \Delta t \to 0 \\ n_t \Delta t = t}}\sum_{n=0}^{n_t-1} \mathbb{E}[ C_{Y_t}^{-1} \Bigg(  \sum_{q=1}^m U_{h \Delta t} b_q( n \Delta t, X_{n \Delta t}) b_q\bigl(n \Delta t,X_{ n \Delta t}\bigr)^\top \Bigg) P_{U_{n \Delta t}}^\perp ] \Delta t, \\
		\end{aligned}
	\end{equation*}
	where in the last line we employ the fact that $U^M_t$ is adapted and \eqref{eq: brownian incerement}.
	From the above relation, we see that (informally) at the limit, $U$ satisfies a deterministic ODE, and hence it is a deterministic quantity. Using this property and the fact that the covariance $C_{Y_t}$ is deterministic, we found that
		\begin{equation*}
		\begin{aligned}
			&\lim\limits_{M \to \infty}\sum_{q=1}^m \sum_{i=1}^{M} \int_{0}^{t}G^U_{i,q} \mathrm{d}W_t^{(i,q)} \\
			=&	 \lim\limits_{\substack{ \Delta t \to 0 \\ n_t \Delta t = t}}\sum_{n=0}^{n_t-1} C_{Y_{n \Delta t}}^{-1} U_{n \Delta t}  \mathbb{E}[\Bigg(  \sum_{q=1}^m b_q(n \Delta t, X_{n \Delta t}) b_q\bigl(n \Delta t,,X_{n \Delta t}\bigr)^\top \Bigg) ]P_{U_{n \Delta t}}^\perp  \Delta t \\
				=&	 \lim\limits_{\substack{ \Delta t \to 0 \\ n_t \Delta t = t}}\sum_{n=0}^{n_t-1} C_{Y_{n \Delta t}}^{-1} U_{n \Delta t}  \mathbb{E}[\Bigg(  \sum_{q=1}^m b_q(n \Delta t, X_{n \Delta t}) b_q\bigl(n \Delta t,X_{n \Delta t}\bigr)^\top \Bigg) ] P_{U_{n \Delta t}}^\perp  \Delta t \\
					=&	 \lim\limits_{\substack{ \Delta t \to 0 \\ n_t \Delta t = t}}\sum_{n=0}^{n_t-1} C_{Y_{n \Delta t}}^{-1} U_{n \Delta t}  \mathbb{E}[b(n \Delta t, X_{n \Delta t}) b\bigl(n \Delta t,X_{n \Delta t}\bigr)^\top] P_{U_{n \Delta t}}^\perp  \Delta t \\
					= & \int_{0}^{t} C_{Y_s}^{-1} U_{s}  \mathbb{E}[b(s, X_{s}) b\bigl(s,X_{s}\bigr)^\top] P_{U_s}^\perp  \mathrm{d}s
		\end{aligned}
	\end{equation*}
	Therefore, in the limit in $M$ the noise term for $U^M$ becomes the following deterministic drift correction for $U$
	\begin{equation*}
	 \sum_{i=1}^M\sum_{q=1}^m \bigl(Y_t\Sigma_M(Y_t^M)^\top\bigr)^{-1}  Y_t\Sigma_MB_{i,q}(t,X_t)^\top P_{U_t}^\perp\mathrm{d}W_t^{(i,q)}  \underset{M\to+\infty}{\longrightarrow}	C^{-1}_{Y_t}\sum_{q=1}^mU_t\,\mathbb{E}\bigl[b_q(t,X_t)b_q(t,X_t)^\top\bigr]P_{U_t}^\perp.
	\end{equation*}
	
	\subsection{The DLRA equations based on Stratonovich derivation}\label{sec: DLRA eq strat}
	We finally obtain the limiting DLRA system in $L^2(\Omega,\mathbb{R}^d)$:
	\begin{equation}\label{eq: new DLRA stra} 
		\begin{aligned}
			\mathrm{d}Y_t &= \left( U_t \Big(a(t,X_t) - \frac{1}{2} \partial_x b(t,X_t) [b(t,X_t)]\Big)+ \frac{1}{2} C_{Y_t}^{-1} \mathbb{E}[Y_t b(t,X_t)^{\top}P_{U_t}^{\perp}] P_{Y_t}^{\perp}[b(t,X_t)] \right) \mathrm{d}t+ U_tb(t,X_t)\mathrm{d}W_t, \\
			\mathrm{d}U_t &= \Big(C_{Y_t}^{-1}\mathbb{E}\bigl[Y_t\Big(a(t,X_t) - \frac{1}{2} \partial_x b(t,X_t) [b(t,X_t)]\Big)^\top\bigr]+ U_t\mathbb{E}\bigl[b(t,X_t)b(t,X_t)^\top\bigr]\\
			& \quad - \frac{1}{2} U_t \mathbb{E}[b(t,X_t) b(t,X_t)^{\top}] + U_t  \mathbb{E}[b(t,X_t)Y_t^{\top}] C_{Y_t} ^{-1} \mathbb{E}[Y_t b(t, X_t)^{\top}]  \Big) P_{U_t}^\perp  \mathrm{d}t, \\
		\end{aligned}
	\end{equation}
	where $X_t = U_t^\top Y_t$ and $C_{Y_t}=\mathbb{E}\bigl[Y_tY_t^\top\bigr]$ is the Gramian of the stochastic basis.
   Finally, the equation for $X$ can be retrieved through Itô formula.
	
	Local and global well-posedness of \eqref{eq: new DLRA stra} can be achieved similarly to results stated in \cite{kazashi2025dynamical} under assumptions of Lipschitzianity and linear-growth bound.

\section*{Conclusion}
\addcontentsline{toc}{section}{Conclusion}
In this article, we introduced alternative DLRA formulations for SDEs, either through minimization of a quantity of interest or by employing the Stratonovich calculus, obtaining different equations to the one proposed in \cite{kazashi2025dynamical}. 

The former approach derives an additional term in the drift resembling a projection onto the horizontal component of the tangent space, i.e.\ on the corange of the surrogate. These equations are equivalent to the one derived in \cite{cao2018stochastic} and can be also generalized to a three-terms DLRA, i.e.\ defined by mean, deterministic, and stochastic bases. 

The latter approach considers a particle approximation of an SDE in a Stratonovich formalism and exploits this setting to express its approximation constrained to low-rank manifold. DLRA equations are retrieved by standard conditions on the basis and Monte-Carlo convergence over the number of samples. The sought limit equations show additional terms in the drift involving the diffusion term and its derivative, restricted to the tangent space onto the surrogate point in the manifold.

These results leave open the question of which DLRA technique is ultimately preferable, setting the next step in this line of research. 

\section*{Acknowledgements}
This work has also been supported by the Swiss National Science Foundation under the
Project n. 200518 “Dynamical low rank methods for uncertainty quantification and data assimilation”.

\printbibliography

\end{document}